\documentclass{article}
\usepackage[english]{babel}
\usepackage[letterpaper,top=2cm,bottom=2cm,left=3cm,right=3cm,marginparwidth=1.75cm]{geometry}
\usepackage{amsmath}
\usepackage{amssymb}
\usepackage{amsthm}
\usepackage{amsfonts}
\usepackage{algorithm}
\usepackage{algorithmicx}
\usepackage{algpseudocode}
\usepackage{authblk}
\usepackage{comment}
\usepackage{enumitem}
\usepackage{float}
\usepackage{bbm}
\usepackage{graphicx}
\usepackage{comment}
\usepackage{todonotes}
\usepackage[colorlinks=true, allcolors=blue]{hyperref}

\allowdisplaybreaks

\newtheorem{theorem}{Theorem}
\newtheorem{definition}{Definition}
\newtheorem{lemma}{Lemma}
\newtheorem{corollary}{Corollary}
\newtheorem{proposition}{Proposition}

\newcommand{\ceil}[1]{{\left\lceil{#1}\right\rceil}}

\newcommand{\epsfin}{\ensuremath{\epsilon_{\textrm{final}}}}
\newcommand{\bigO}[1]{\mathcal{O}\left( #1 \right)}

\newcommand{\bigTh}[1]{\Theta\left( #1 \right)}
\newcommand{\bigOm}[1]{\Omega\left( #1 \right)}
\newcommand{\zeroes}{\ensuremath{\mathbf{0}}}
\newcommand{\ones}{\ensuremath{\mathbf{1}}}

\title{Solution Size Growth as an Obstacle to Iterative Refinement for Linear Optimization}

\author[1]{Henry Graven\thanks{hgraven@usc.edu}}
\author[1,2]{Giacomo Nannicini\thanks{g.nannicini@usc.edu}}

\affil[1]{Daniel J.~Epstein Department of Industrial and Systems Engineering, University of Southern California}
\affil[2]{Ming Hsieh Department of Electrical and Computer Engineering, University of Southern California}

\begin{document}
\maketitle

\begin{abstract}
  Iterative refinement is a technique used to exponentially boost the precision of the solution to certain problems by solving a sequence of related problem instances. We study the size of the solutions arising from the sequence of problems produced by the iterative refinement algorithm when applied to linear optimization problems (LPs). We introduce a class of problems for which the size of such solutions can grow exponentially with the number of refinement iterations, and we show evidence that to avoid this undesirable behavior, a black-box LP solver would have to be \emph{lucky}---a concept introduced formally in this paper, by which we mean that the solver would have to return a solution of better quality than what is required by its termination tolerances. Furthermore, we show that the exponential scaling occurs if one uses a natural instantiation of the simplex algorithm using Bland's pivoting rule. We discuss the implications of our lower bound on the solution size for optimization algorithms whose running time depends polynomially on the solution size, such as quantum multiplicative weights update algorithms and quantum matrix games algorithms. 
\end{abstract}

\section{Introduction}
\label{sec:intro}

The time complexity of algorithms to solve continuous optimization problems typically depends on the target tolerance $\epsilon$ of the solution; the tolerance may refer to a feasibility tolerance, an optimality tolerance, or both. This paper is concerned with linear programs (LPs---as is customary, we use the acronym LP to also denote ``linear programming''); for ease of exposition, we use the same $\epsilon$ to refer to primal and dual feasibility, as well as the complementarity gap (to be defined later). \emph{Iterative refinement} (we use IR as shorthand) is a technique that can be used to boost the precision of a solution by solving a polylogarithmic number of related problem instances. It was first developed in the context of solving linear systems \cite{wilkinson1963rounding,wilkinson1994rounding}, where the idea is easily explained: suppose we solve a linear system to fixed precision, say, two digits after the decimal dot. After obtaining a solution with those two digits correct, we compute the \emph{residual} to assess the error of the current solution with respect to an exact solution of the linear system.  With the residual, we translate and possibly rescale the problem, obtaining a new linear system (called the \emph{refining problem}) whose solution is a correction term for the initial solution that gives an additional two digits after the decimal dot. This brings the total to four digits, once added to the initial solution. The process can be iterated. If the initial solution satisfies the error tolerance $\epsilon$, after $k$ iterations we obtain a solution that satisfies error tolerance $\epsilon^k$. Thus, only $\ceil{\log_{\epsilon} \epsfin}$ iterations of IR are necessary to satisfy tolerance $\epsfin$. This idea can be extended to linear programs \cite{gleixner2012improving,gleixner2016iterative,gleixner2020linear} and to semi-definite programs \cite{mohammadisiahroudi2025quantum}, exploiting linearity of the constraints. In the context of linear and semi-definite optimization, each refining problem is an optimization problem, and it gives a correction to an existing solution that decreases its primal and dual infeasibilities, and the complementarity gap.

This paper studies the \emph{size} of the solutions of the optimization problems that arise in the course of the application of IR for LP. Our assumption, similar to that of \cite{gleixner2020linear}, is that we want to solve a given LP with small error tolerance $\epsfin \approx \epsilon^k$, using a \emph{limited-precision LP oracle}: an algorithm to solve LPs (often called an \emph{LP solver}) that outputs a solution that is only guaranteed to be approximately feasible and approximately optimal, i.e., with tolerance $\epsilon$. Here and in the following, $\epsfin$ refers to the target error tolerance to be achieved after running the IR algorithm, whereas $\epsilon$ refers to the tolerance of the available limited-precision LP oracle. IR can be directly applied to this setting. In the standard IR algorithm for LP \cite{gleixner2016iterative,gleixner2020linear}, the bounds for the variables of the refinement LP solved in every iteration may grow exponentially in the number of iterations. This yields an exponential upper bound to the size $R^{(t)}$ (defined as the largest entry of the solution vector) of the solution of the refining problem at iteration $t$, $t=1,\dots,k$. This paper is concerned with providing a lower bound to $R^{(t)}$: do there exist LPs for which the sequence of solutions encountered when applying IR has exponentially-growing size in $t$? Here and in the following, we use a superscript with parentheses, e.g., ${R}^{(t)}$ to denote a quantity at some iteration $t$, whereas exponential functions do not have brackets, e.g., ${R}^{t}$.

In the context of commonly-used algorithms for linear programming, such as the simplex algorithm or interior point methods, exponentially-growing solution sizes would have almost negligible effect. This is because the running time of those algorithms depends polylogarithmically on the solution size, so even exponentially-large solutions can be accommodated at little extra cost: if the solution size doubles in every IR iteration, we need an extra bit of precision per iteration (if the solution size grew at a doubly-exponential rate, that would create difficulties). The total number of iterations of the IR algorithm is logarithmic in $\epsfin$, so even (single) exponential growth of $R^{(t)}$ can be accommodated by using $\bigO{\log(1/\epsfin)}$ additional bits. However, in recent years the study of quantum optimization algorithms has produced multiple examples of quantum algorithms that can solve optimization problems with a running time that has attractive dependence on the dimension of the problem, but polynomial dependence on the inverse tolerance $1/\epsilon$ as well as the largest entry $R$ of the solution vector (and, possibly, a similar parameter related to the size of the dual solution). Examples of such algorithms are \cite{vanApeldoorn2019,apeldoorn2017quantum,apeldoorn2020sdp,brandao2017quantum,brandao2019quantum,li2020sublinear,augustino2023central}; several of these have a running time that depends sublinearly on the dimension of the problem, but $R$ and $1/\epsilon$ appear in the running time expression with a polynomial of possibly high degree. In addition to the quantum algorithms mentioned above, there are classical algorithms that exhibit similar dependence in their asymptotic running time expression. Examples are classical analogues of quantum algorithms, such as the multiplicative weights update algorithm \cite{arora2012multiplicative} and matrix games algorithms \cite{arora2012multiplicative,vanApeldoorn2019}, as well as first-order algorithms for saddle point problems \cite{nemirovski2004prox} when applied to LPs.

In fact, $R$ and $\epsilon$ can be traded against each other by rescaling the problem, so to obtain an unscaled solution with tolerance $\epsilon$ we can equivalently scale the problem data down by a factor $R$, set the error tolerance to $\epsilon/R$, and the rescaled problem has $R=1$. This type of running time dependence on $R$ often occurs if there is some relative error bound (e.g., optimality gap $\epsilon_r \|x^{\star}\|$, where $x^{\star}$ is the optimal solution) that we want to convert to an absolute error bound (e.g., optimality gap $\epsilon_a$), which is typically used in the optimization world. Note that because of this, we are interested in the actual solution size (i.e., some norm of the solution, usually the infinity norm) rather than the length of its bit encoding: a small solution, say, of size $c/2^k$ for some number $c < 2^k$, is typically encoded using $\bigO{k}$ bits, and to achieve absolute error $\epsilon_a$ we choose $\epsilon_r = \bigO{\epsilon_a}$; on the other hand, a large solution, say, of size $2^k$, also has bit encoding of length $\bigO{k}$, but to achieve absolute error $\epsilon_a$ we choose $\epsilon_r = \bigO{\epsilon_a/2^k}$, which is equivalent to rescaling the problem so that the optimal solution has unit norm. 

It is natural to think of IR as a possible approach to improve the running time dependence of the algorithms discussed above on the tolerance parameter $\epsfin$, from polynomial to polylogarithmic in $1/\epsfin$. If we could apply IR and guarantee that the size parameter $R^{(t)}$ remains approximately constant $= \bigO{R}$ through the iterations, then we could reach a solution with final tolerance $\epsfin$ after $\bigO{\log(1/\epsfin)}$ iterations of the IR algorithm, each executed with $1/\epsilon = \bigO{R}$. If $R^{(t)}$ grows polynomially in $t$, the use of IR may still yield benefits. However, if $R^{(t)}$ grows exponentially in $t$, then we need to choose $1/\epsilon = \bigO{R^{(t)}} = \bigO{R^t}$ at iteration $t$. Because $t = \bigO{\log(1/\epsfin)}$ iterations, this means that the final iterations of the IR scheme are executed with $\epsilon = \bigO{\epsfin}$: we do not achieve the exponential advantage in precision that IR promises.

The main result of the paper is to show that in pathological cases the solution size of IR problems can grow exponentially with the number of iterations. We take the LP solver as a black box that delivers a basic solution with error tolerance $\epsilon$, where $\epsilon$ can be specified. Because we do not insist on the specifics of the LP solver, we do not give an unconditional proof of exponential scaling; however, we provide evidence that the exponential scaling appears difficult to get around. Specifically, we construct a family of problems, parameterized by $\epsilon$ and an integer $T$, with the following property: starting from a vertex that is primal feasible and with dual infeasibility $\bigTh{\epsilon^t}$, a solution of the refining problem the leads to a vertex with dual infeasibility $\bigTh{\epsilon^{t+1}}$ has size an order of magnitude larger than the previous one. As a consequence, either the LP solver gets \emph{lucky}, i.e., it finds a solution with dual infeasibility much smaller than it can guarantee given its tolerance $\epsilon$, or it returns solutions that get exponentially larger as the iterations progress. We also show that if one applies the simplex algorithm with Bland's pivoting rule to the refining problems obtained from this family of LPs, the sequence of refining solutions incurs the pathological case.

As a consequence of our results, IR cannot be applied ``as is'' to improve the running time of the classical and quantum algorithms referenced above in terms of the tolerance parameter $\epsilon$; or, more in general, of any algorithm whose running time depends polynomially on the largest entry of the solution. This is particularly important for the prospects of sublinear-time quantum algorithms, because it indicates that we cannot retain sublinear time in the problem dimension while achieving polylogarithmic $\epsilon$ dependence thanks to IR.

Our results hold under certain conditions. To show that a solver either gets lucky or incurs into the pathological case, we assume that the solver returns basic solutions and that degeneracies in the bases are resolved in a certain way (i.e., if a vertex is degenerate, the solver returns a corresponding basis that favors a specific set of variables to be in the basis in degenerate rows; however, as the problem size increases, we show that the fraction of degenerate vertices goes to zero). For the pathological example using Bland's rule, we choose a natural starting basis for the simplex algorithm applied to the refining problems. Thus, the exponential lower bound on the solution sizes shown here does not hold for every LP solver; for example, an interior point method may not return a basic solution. However, we argue that our results are already difficult to get around: an interior point \emph{could} return a basic solution, and therefore the sequence of solutions presented in this paper. In other words, our goal is to provide a hardness result that makes it difficult to show that an LP solver \emph{always} avoids the pathological example: as long as the solver may sometimes return a basic solution, the exponential growth of the IR solutions shown in this paper implies that IR does not always reduce the $\epsilon$ dependency to polylogarithmic. Moreoever, it is plausible that our family of LPs can be modified to produce pathological examples tailored to a given LP solver. This situation reminds us of the proofs for the exponential running time of the simplex algorithm: Klee and Minty's famous example in the 70s \cite{klee1972good} gave a lower bound that only applied to Dantzig's pivoting rule, and many other examples for different pivoting rules quickly followed, e.g., \cite{avis1978notes,goldfarb1979worst,goldfarb1994complexity,jeroslow1973simplex}. More than 50 years later, there is no single proof of the exponential running time of the simplex algorithm that applies to all its variants, but this topic is still an active area of research, see, e.g., \cite{disser2023exponential,black2024exponential}.

\paragraph{Overview of the paper.}
In Section~\ref{sec:preliminaries} we introduce our notation and give a formal description of the IR algorithm. In Section~\ref{sec:results} we give a formal description of the LP family, and we state our main results. Subsequent sections contain a detailed proof of the statements in Section~\ref{sec:results}. More specifically, Section~\ref{sec:vertexstruct} discusses the vertex structure of the polytope, Section~\ref{sec:lucky-convergence} proves our first main theorem, Section~\ref{sec:simplex} proves our second main theorem. In the remainder of this section we give an informal description of our results; the description is not indended to be accurate, but it provides a good intuition of the structure of the proofs.

\paragraph{Informal summary of the proof.}
Fix an integer $T$. We describe a problem that has two \emph{controller} variables $w_1, w_2$, and $T+1$ \emph{structural} variables $z_1,\dots,z_{T+1}$. All variables are nonnegative and there are two simplex-like constraints: the $w$ variables sum up to $1$, and the $z$ variables sum up to $1$. The remaining constraints enforce the following vertex structure: at each vertex, either odd-index or even-index $z$ variables are nonzero. We have odd-index vertices whenever $w_1 \in \{1, \frac{1}{3}, \frac{1}{6}, \dots, \frac{1}{3(T+1)}\}$, and even-index vertices whenever $w_2 \in \{1, \frac{1}{3}, \frac{1}{6}, \dots, \frac{1}{3(T+1)}\}$. When $w_1 \in \{1, \frac{1}{3}\}$ the odd-index $z$ variables form a standard simplex, and similarly when $w_2 \in \{1, \frac{1}{3}\}$ for the even-index $z$ variables. At remaining values of the controller variables, the vertices split a unit mass among multiple odd- or even-index $z$ variables. The $z$ variables and their objective function coefficients are scaled in such a way that each variable is slightly more attractive than the previous one in terms of objective function contribution. Crucially, there is an improving path on the edges of the polytope that starts at the vertex with $z_1=1$, then proceeds to the vertex with $z_2=1$, and so on until reaching the vertex with $z_{T+1}=1$. (In fact, the $z$ variables are scaled and their value at the vertices is not $1$; in this simplified overview we ignore this aspect.) To follow this path, the values of the controller variables $w$ need to ``bounce'' because we alternate between nonzero odd- and even-index $z$ variables; such a bounce changes the $w$ variables by at least a constant (for $T$ fixed). If the path is discovered by solving the refining problem, the bounce requires the $w$ variables to get larger and larger, due to the fact that the refining problem zooms in by a factor $\epsilon$ in each iteration. This causes the exponential scaling. By construction, Bland's pivoting rule follows the aforementioned path on the polytope, which shows one of our results. The second result is more involved: we show that in each iteration of IR, obtaining a solution with a dual infeasibility that is ``just right'' (not too small and not too large compared to the tolerance $\epsilon$) also requires alternating between nonzero odd- and even-index $z$ variables, which in turn requires the $w$ variables to bounce. This leads to the conclusion that, unless the LP solver gets lucky and determines a solution of higher quality than what it is tasked to do, the IR algorithm produces refining solution with exponentially-growing size.

\section{Preliminaries}
\label{sec:preliminaries}
We begin our study by introducing our notation and conventions. We use $\|x\|, \|x\|_1, \|x\|_\infty$ to denote the $\ell^2$, $\ell^1$, $\ell^\infty$ norm of a vector, i.e., the Euclidean norm, the sum of the absolute value of the entries, and the absolute value of the largest entry. We denote by $\zeroes$ the all-zero vector and by $\ones$ the all-ones vector, of size determined by the context. We often index matrices and vectors using the name of a decision variable rather than its position, which simplifies the exposition especially when dealing with bases of an LP, so that we do not have to keep track of the order of the elements in the basis; e.g., if we have a set $B = \{w_1, w_2, z_1, z_3\}$ and a vector $\lambda$ of size $|B|$, we write $\lambda_{z_1}$ to refer to the element of $\lambda$ corresponding to $z_1$. We use this notation only when it does not create ambiguities.

The LP family used in this paper depends on an integer $T$ and a rational number $\epsilon > 0$. Unless stated otherwise, $T$ is fixed and the hidden constants may depend polynomially on $T$, but not on $\epsilon$. That is, all asymptotic statements involving $\epsilon$ in the main theorem are interpreted with $T$ fixed, and in particular they may only hold for $\epsilon$ small enough where the threshold for ``small enough'' depends on $T$ and other constants. Thus, for positive functions $f(\epsilon)$ and $g(\epsilon)$, we write $f=\bigO{g}$ if there exist constants $\mu_T>0$ and $\epsilon_0>0$, where $\mu_T$ may depend polynomially on $T$, such that $f(\epsilon) \le \mu_T g(\epsilon)$ for all $0 < \epsilon < \epsilon_0$. Similarly, $f=\bigOm{g}$ means that there exist constants $\beta_T>0$ and $\epsilon_0>0$, where $\beta_T$ may depend polynomially on $T$, such that $f(\epsilon) \ge \beta_T g(\epsilon)$ for all $0 < \epsilon < \epsilon_0$. $f=\bigTh{g}$ means both $f=\bigO{g}$ and $f=\bigOm{g}$.

It is important to remark why the assumption that $T$ is fixed does not weaken our conclusions. Typically, in an application of the IR algorithm the input parameters are the refining tolerance $\epsilon$ and the final tolerance $\epsfin \ll \epsilon$. To solve the problem in $T$ IR iterations we need $\epsilon \le \sqrt[T]\epsfin$, so smaller $\epsfin$ implies smaller $\epsilon$. For every integer $t=1,\dots,T$, our results show that there is some value of $\epsilon$ below which the size of the solutions to the refining problem is $\bigOm{1/\epsilon^t}$. Because this holds for every $T$, this implies that we cannot achieve arbitrarily small error $\epsfin$ without exponential growth of the solutions. We also note that $T = \bigO{\text{poly} \log(1/\epsfin)}$, so if $T$ is not considered fixed our results still hold up to polylogarithmic factors.

To introduce the IR algorithm we first need to define the LP problem to which it is applied and its dual. For the standard-form LP
\begin{equation}
  \label{eq:genericplp}
\min\{c^\top x: Ax = b,\ x \ge \ell\}
\end{equation}
the dual is
\begin{equation*}
\max\{b^\top y + \ell^\top s(y): A^\top y + s(y) = c,\ s(y) \ge \zeroes \};
\end{equation*}
since the dual slack vector $s(y)$ is uniquely determined by the dual variables $y$, we denote it as a function of $y$, i.e., $s(y) := c - A^\top y$. (This also helps to distinguish from the slack vector $s$ used after Section~\ref{sec:results}.) We use
\begin{equation*}
\gamma(x,y):=(x - \ell)^\top s(y) = (x - \ell)^\top(c - A^\top y)
\end{equation*}
for the complementary-slackness residual; this is positive for feasible primal-dual pairs, but it could be negative in case of infeasibilities. The corresponding complementarity violation is $|\gamma(x,y)|$. We say that a (not necessarily feasible) primal-dual pair is complementary slack if $|\gamma(x,y)| = 0$; for feasible primal-dual pairs this is equivalent to the usual component-wise condition
\begin{equation*}
(x_j - \ell_j)(s(y))_j=0 \qquad \forall \ j.
\end{equation*}

Among all solutions that solve a given LP, basic solutions play an important role in theory and in practice.
\begin{definition}[Basic primal-dual solution]
\label{def:basic-primal-dual}
Let $A\in\mathbb{Q}^{m\times n}$ be of full row rank. A primal-dual pair $(x,y)$ for \eqref{eq:genericplp} is called \emph{basic} if there exists a set of column indices $B\subseteq\{1,\dots,n\}$ with $|B|=m$ such that the submatrix $A_B$ is non-singular,
\begin{equation*}
x_j = \ell_j \qquad \forall \ j \notin B,
\end{equation*}
and $x_B$ and $y$ are determined by
\begin{equation*}
A_B x_B = b - \sum_{j \notin B}A_j\ell_j, \qquad y^\top A_B = c_B^\top.
\end{equation*}
\end{definition}
A basic primal-dual pair need not be primal feasible. We write \emph{primal-feasible basic pair} for a basic pair such that its primal component satisfies the constraints of \eqref{eq:genericplp}. A primal-feasible basic pair corresponds to a vertex of the feasible region of \eqref{eq:genericplp}.

In this paper we consider a slightly modified version of the IR algorithm presented in \cite{gleixner2020linear}, formally stated in Algorithm~\ref{alg:ir}. In our version, the residual and scaled data at every iteration of the algorithm are formed exactly; in \cite{gleixner2020linear}, they are allowed to be computed in finite precision. Note that a solution to the refining problem formed with exact data is an approximate solution to its finite-precision representation, with a relative error $\bigO{\epsilon}$ if we truncate the data with error $\epsilon$; thus, we believe that our results also apply to the finite-precision model of \cite{gleixner2020linear} with some modifications, although we did not pursue a detailed study of this direction. In practice, it is reasonable to assume that the data of a problem can be stored, and calculations can be carried out, with much smaller error than the error tolerance $\epsilon$ of the LP solver; this is typically what happens in virtually every existing LP solver, as a floating point representation already provides much higher precision than the feasibility or optimality tolerance used by the solver.

\begin{algorithm}[ht]
\caption{Iterative refinement with residuals computed exactly} \label{alg:ir}
\begin{algorithmic}[1]
\Require $\min \{ c^{\top} x \mid Ax = b, x \geq \ell \}$ with $A \in \mathbb{Q}^{m \times n}$, $b \in \mathbb{Q}^{m}$, $c, \ell \in \mathbb{Q}^{n}$, termination tolerance $\epsfin > 0$, refining tolerance $\epsilon > \epsfin$, incremental scaling limit $\alpha$
\Ensure Primal-dual solution $x^* \in \mathbb{Q}^{n}$, $y^* \in \mathbb{Q}^{m}$
\State $\Delta_P^{(0)} \gets 1,\ \Delta_D^{(0)} \gets 1$
\State Set $\hat{b} \leftarrow b, \hat{\ell} \leftarrow \ell, \hat{c} \leftarrow c$ \label{line:getinitdata} 
\State Solve $\min \{ \hat{c}^{\top} x \mid A x = \hat{b}, x \geq \hat{\ell} \}$ with tolerance $\epsilon$ \label{line:cumulprob} \Comment{initial solve}
\State $x^{(0)},y^{(0)} \gets$ $\epsilon$-approximate primal-dual solution returned
\For{$t = 1, 2, \dots$} \Comment{refinement loop}
    \State $b^{\text{res}} \gets b - A x^{(t-1)}$ \Comment{compute violations}
    \State $\ell^{\text{res}} \gets \ell - x^{(t-1)}$
    \State $c^{\text{res}} \gets c - A^{\top} y^{(t-1)}$
    \State $\delta_{P}^{(t-1)} \gets \max \{ \max_{j=1,...,m} | b^{\text{res}}_{j} |, \max_{i=1,...,n} \ell^{\text{res}}_{i} \}$ \Comment{check termination}
    \State $\delta_{D}^{(t-1)} \gets \max \{ 0, \max \{ -c^{\text{res}}_{i} \mid i = 1, \dots, n \} \}$
    \State $\delta_{S}^{(t-1)} \gets \| \sum_{i=1,...,n} \ell^{\text{res}}_{i} \cdot c^{\text{res}}_{i} \|$
    \If{$\delta_{P}^{(t-1)}\le \epsfin$, $\delta_{D}^{(t-1)}\le \epsfin$, and $\delta_{S}^{(t-1)} \leq \epsfin$}
        \State \Return $x^* \gets x^{(t-1)}$, $y^* \gets y^{(t-1)}$
    \EndIf
    \State \label{line:deltapdef} $\Delta_{P}^{(t)} \gets 1/\max\{ \delta_{P}^{(t-1)}, (\alpha \Delta_{P}^{(t-1)})^{-1} \}$
    \State \label{line:deltaddef} $\Delta_{D}^{(t)} \gets 1/\max\{ \delta_{D}^{(t-1)}, (\alpha \Delta_{D}^{(t-1)})^{-1} \}$
    \State \label{line:gethats} $\hat b \gets \Delta_{P}^{(t)}b^{\mathrm{res}},\quad \hat\ell \gets \Delta_{P}^{(t)}\ell^{\mathrm{res}},\quad \hat c \gets \Delta_{D}^{(t)}c^{\mathrm{res}}$
    \State Solve $\min \{ \hat{c}^{\top} x \mid A x = \hat{b}, x \geq \hat{\ell} \}$ with tolerance $\epsilon$ \Comment{solve transformed problem} \label{line:refprob}
    \State $\hat{x}^{(t)}, \hat{y}^{(t)} \gets$ $\epsilon$-approximate primal-dual solution returned 
    \State \label{line:updatesol} $x^{(t)} \gets x^{(t-1)} + \frac{1}{\Delta_{P}^{(t)}}\hat{x}^{(t)} \quad \text{and} \quad$ $y^{(t)} \gets y^{(t-1)} + \frac{1}{\Delta_{D}^{(t)}} \hat{y}^{(t)}$ 
\EndFor
\end{algorithmic}
\end{algorithm}

%% We work with the exact-data form of IR: the residuals, scaling factors, and transformed LP data are computed and passed to the oracle exactly. This deliberate and reasonable restriction of the general finite-precision implementation considered in \cite{gleixner2016iterative} is outlined in the assumption below.

%% \begin{assumption}[Exact-data primal-feasible basis oracle]
%% \label{ass:exact-data-oracle}
%% Algorithm~\ref{alg:ir}
%% is interpreted in an exact-data model. At every refinement iteration
%% $t$, the residuals and the scaled data are formed exactly:
%% \begin{equation*}
%% \hat b^{(t)}=\epsilon^{-t}(b-Ax^{(t-1)}),\qquad
%% \hat\ell^{(t)}=\epsilon^{-t}(\ell-x^{(t-1)}),\qquad
%% \hat c^{(t)}=\epsilon^{-t}(c-A^\top y^{(t-1)}).
%% \end{equation*}
%% The oracle returns a basis together with a primal-dual pair that is
%% basic for this exact refining LP, whose primal component is exactly
%% feasible for that refining LP, and whose dual component satisfies the
%% dual feasibility tolerance
%% \begin{equation*}
%% \hat c^{(t)}-A^\top \hat y^{(t)}\ge -\epsilon\mathbf 1.
%% \end{equation*}
%% The initial solve satisfies the analogous conditions for the original
%% LP.
%% \end{assumption}

On lines \ref{line:deltapdef} and \ref{line:deltaddef} we define scaling factors $\Delta_{P}^{(t)}, \Delta_{D}^{(t)}$ for the primal and dual problem computed from the residual vectors. If we choose $\alpha = 1/\epsilon$ as the incremental scaling limit, then it is easily seen by induction that $\Delta_{P}^{(t)} = \Delta_{D}^{(t)} = 1/\epsilon^t$ at every iteration. We make such a choice in the rest of this paper, without loss of generality, as it simplifies the exposition.

\begin{theorem}[Corollary 3.6 of \cite{gleixner2016iterative}]
  \label{thm:irconvergence}
  Let $\epsilon \in [0,1)$. Assume that we have access to an LP solver that, given problem \eqref{eq:genericplp}, returns a primal-dual pair $(x,y)$ such that
  \begin{align*}
    \|Ax -b\|_\infty &\le \epsilon, \\
    x-\ell &\ge -\epsilon \ones, \\
    c-A^\top y &\ge -\epsilon \ones, \\
    |\gamma(x,y)| &\le \sigma.
  \end{align*}
  Then, for every refinement iteration $t$, the primal-dual iterate
$(x^{(t)},y^{(t)})$ computed by Algorithm~\ref{alg:ir} satisfies
  \begin{align*}
    \|Ax^{(t)}-b\|_\infty &\le \epsilon^{t+1}, \\
    x^{(t)}-\ell &\ge -\epsilon^{t+1} \ones, \\
    c-A^\top y^{(t)} &\ge -\epsilon^{t+1} \ones, \\
    |\gamma(x^{(t)},y^{(t)})| &\le \sigma \epsilon^{2t}.
  \end{align*}
  In particular, if we choose $\epsilon = \sigma$ to be a constant $< 1$, Algorithm~\ref{alg:ir} reaches any prescribed final tolerance $\epsfin$ after $\bigO{\log(1/\epsfin)}$ refinement iterations.
\end{theorem}
In the following we need a measure of dual infeasibility for the iterates $x^{(t)}, y^{(t)})$ computed by Algorithm~\ref{alg:ir}. For $t=0,1,\dots$, we define it as
\begin{equation*}
\delta_D^{(t)} := \left\| \min\left\{ \zeroes,c-A^\top y^{(t)} \right\} \right\|_\infty.
\end{equation*}
A similar measure for the refining problem is also helpful; for
$t=1,2,\dots$ we define it as
\begin{equation*}
\widehat \delta_D^{(t)} := \left\| \min\left\{ \zeroes, \hat c^{(t)}-A^\top \hat y^{(t)} \right\} \right\|_\infty.
\end{equation*}
Note that these quantities are unambiguously defined by the sequences $(x^{(t)}, y^{(t)})$, $(\hat{x}^{(t)}, \hat{y}^{(t)})$ computed by Algorithm~\ref{alg:ir}.

It is important to distinguish between two types of problems involved in the IR algorithm: the \emph{cumulative} problem, i.e., the original problem \eqref{eq:genericplp} that we wish to solve, and for which an initial solution is computed on line \ref{line:cumulprob} of Algorithm~\ref{alg:ir}; and the \emph{refining} problems, obtained from the scaled residuals in each iteration and solved on line \ref{line:refprob} of Algorithm~\ref{alg:ir}. We study IR through the cumulative problem, but analyze its effects on the associated refining problems and their solutions. This is possible because the cumulative and refining problems are linked by affine transformations: a refining solution at iteration $t$ corresponds to an updated cumulative iterate, with the relevant feasibility and optimality errors scaled by powers of $\epsilon$. Thus, we can formulate statements directly for the cumulative problem and interpret their implications on the refining problems. Recall that $\Delta_{P}^{(t)} = \Delta_{D}^{(t)} = 1/\epsilon^t$ for every $t$. The correspondence between the cumulative and refining problems is characterized in the next proposition. Note that we generally use $\hat{\cdot}$ to denote quantities related to the refining problem. 
\begin{proposition}[Correspondence between the cumulative and refining problems]
\label{prop:cumulative-refining-correspondence}
Let $(x^{(t-1)},y^{(t-1)})$ be the primal-dual solution for the cumulative problem computed by Algorithm~\ref{alg:ir} at iteration $t-1$. Let
\begin{equation*}
\hat b^{(t)} = \epsilon^{-t} \left( b-Ax^{(t-1)} \right),
\qquad
\hat \ell^{(t)} = \epsilon^{-t} \left( \ell-x^{(t-1)} \right),
\qquad
\hat c^{(t)} = \epsilon^{-t}\left( c-A^\top y^{(t-1)}  \right),
\end{equation*}
be the corresponding data for the refining problem solved at iteration $t$.  Let $(\hat x^{(t)},\hat y^{(t)})$ be a primal-dual pair for the refining problem $\min \{ \hat{c}^{\top} x \mid A x = \hat{b}, x \geq \hat{\ell} \}$, and let
\begin{equation*}
x^{(t)} := x^{(t-1)} + \epsilon^t\hat x^{(t)}, \qquad
y^{(t)} := y^{(t-1)} + \epsilon^t\hat y^{(t)}.
\end{equation*}
Then the following hold.
\begin{enumerate}
%% \item If $\hat x^{(t)}$ is exactly primal feasible for the refining
%%   problem, then $x^{(t)}$ is exactly primal feasible for the
%%   cumulative problem.
\item If $(\hat x^{(t)},\hat y^{(t)})$ is a basic primal-dual pair for
  the refining problem with basis $B$, then $(x^{(t)},y^{(t)})$ is a
  basic primal-dual pair for the cumulative problem with the same
  basis $B$.
\item $c-A^\top y^{(t)} = \epsilon^t \left( \hat c^{(t)}-A^\top \hat
  y^{(t)} \right)$.
\item If $\hat c^{(t)}-A^\top \hat y^{(t)} \geq -\epsilon\mathbf 1$,
  then $c-A^\top y^{(t)} \geq -\epsilon^{t+1}\mathbf 1$.
\item $\delta_D^{(t)} = \epsilon^t\widehat \delta_D^{(t)}$.
\end{enumerate}
\end{proposition}
\begin{proof}
  Part 1 follows directly from the affine transformation identities of Theorem~3.1 in \cite{gleixner2016iterative}, given our choice of $\Delta_{P}^{(t)}, \Delta_{D}^{(t)}$.

  For part 2, use the definition of the objective for the refining problem:
  $\hat c^{(t)} = \epsilon^{-t} \left(c-A^\top y^{(t-1)} \right).$
  Then
  \begin{equation*}
      \epsilon^t \left( \hat c^{(t)}-A^\top \hat y^{(t)} \right)
      = c-A^\top y^{(t-1)} - \epsilon^t A^\top \hat y^{(t)} 
      = c-A^\top \left( y^{(t-1)}+\epsilon^t\hat y^{(t)} \right) 
      = c-A^\top y^{(t)}.
  \end{equation*}
  Part 3 is a restatement of the same identity. Taking componentwise negative parts and then the infinity norm gives part 4. 
\end{proof}
%% Proposition~\ref{prop:cumulative-refining-correspondence} allows us to move freely between the refining and cumulative viewpoints. In particular, once a property of the refining problem has been established at iteration $t$, it immediately yields a corresponding statement about the cumulative iterate after the $t$-th refinement step.

\section{Description of the LP family and main results}
\label{sec:results}

In this section we formally state our main results, preceded by a description of the family of LPs used for the proofs. Each LP in the family is parametrized by an integer $T\ge 3$ and $\epsilon\in\mathbb{Q}\cap(0,1)$. We use three sets of variables for the primal: controller variables $w$, structural variables $z$, slack variables $s$. The vector $x = (w, z, s)$ collectively denotes all primal variables. The dual variables are still called $y$. We no longer use $s(y)$ for the dual slacks.

%% Note that in Section~\ref{sec:preliminaries} we use $x$ to denote the decision variables of the primal LP, $y$ for the decision variables of the dual LP, and $s(y)$ for the corresponding dual slacks. , but we no longer reserve a symbol for the dual slacks. We find that this notational choice makes the exposition clearer, and because we no longer need to reference the general LP \eqref{eq:genericplp} and its dual, no ambiguity arises. In the few occasions where we still need to refer to the primal-dual sequence $(x^{(t)}, y^{(t)})$ of the IR algorithm, we explicitly refer to Algorithm~\ref{alg:ir}.

\paragraph{Cumulative primal LP.}
%% Let
%% \begin{equation*}
%% z := (x_1, x_2, x_3, \ldots, x_{T+2}, x_{T+3},
%% s_3, \ldots, s_{T+2}, s_{T+3})^\top
%% \in \mathbb{R}^{2T+4}.
%% \end{equation*}
The main problem studied in this paper is
\begin{equation}
\label{eq:PT}
\tag{$P_T$}
\left.
\begin{array}{rrll}
\min
& \displaystyle
\sum_{j=1}^{T}
\bigl(\epsilon^T+\epsilon^{-j}\bigr)z_j +
\bigl(\epsilon^T+\epsilon^{-T}\bigr)z_{T+1} \\
\textnormal{s.t.:} & w_1+w_2 &=& 1, \\
& \displaystyle \sum_{j=1}^{T}\epsilon^{-j}z_j +\epsilon^{-T}z_{T+1} &=& 1, \\
& -3w_1+\epsilon^{-j}z_j+s_j &=& 0 \qquad
\forall j\in\{1,\ldots,T\}\textnormal{ odd}, \\
& -3w_2+\epsilon^{-j}z_j+s_j &=& 0 \qquad \forall j\in\{1,\ldots,T\}\textnormal{ even}, \\
& -3w_1+\epsilon^{-T}z_{T+1}+s_{T+1} &=& 0 \qquad \textnormal{if }T+1\textnormal{ is odd}, \\
& -3w_2+\epsilon^{-T}z_{T+1}+s_{T+1} &=& 0 \qquad \textnormal{if }T+1\textnormal{ is even}, \\
& w, z, s &\ge&\zeroes.
\end{array}
\right\}
\end{equation}
Calling $n$ the number of variables and $m$ the number of constraints
of \eqref{eq:PT}, we have $n = 2T+4$, $m = T+3$. As is customary, we denote the constraint matrix of \eqref{eq:PT} by $A$, the objective function vector by $c$, and the r.h.s.\ vector by $b$.

Temporarily dropping the slack variables $s$ for ease of exposition, the constraints become
\begin{equation*}
\forall j\in\{1,\ldots,T\} \qquad \epsilon^{-j}z_j \leq 3w_1 \textnormal{ if $j$ is odd}, \qquad \epsilon^{-j}z_j \leq 3w_2 \textnormal{ if $j$ is even,}
\end{equation*}
together with
\begin{equation*}
\epsilon^{-T}z_{T+1}\leq 3w_1 \textnormal{ if $T+1$ is odd},
\qquad \epsilon^{-T}z_{T+1}\leq 3w_2 \textnormal{ if $T+1$ is even}.
\end{equation*}
The problem has $T+3$ variables and $T+3$ constraints in this inequality form.

%% Equivalently, in the variable ordering used above, the objective
%% coefficients are
%% \begin{equation*}
%% c_{x_i}
%% =
%% \begin{cases}
%% 0, & i=1,2,\\
%% \epsilon^T+\epsilon^{2-i}, & i=3,\ldots,T+2,\\
%% \epsilon^T+\epsilon^{-T}, & i=T+3,
%% \end{cases}
%% \end{equation*}
%% \begin{equation*}
%% c_{s_j}=0
%% \qquad
%% \forall j\in\{3,\ldots,T+3\},
%% \end{equation*}
%% and
%% \begin{equation*}
%% b=
%% \begin{bmatrix}
%% 1\\
%% 1\\
%% 0\\
%% \vdots\\
%% 0
%% \end{bmatrix}
%% \in\mathbb{R}^{T+3},
%% \qquad
%% \ell=0.
%% \end{equation*}

%% We use $\hat{x}$ to denote decision variables in the refining problem,
%% and $\tilde{x}$ for decision variables in a transformed version of the
%% refining problem where all variable lower bounds are zero; these are
%% properly defined in the Appendix. We use $\hat{c},\hat{l},\hat{b}$ to
%% denote the objective, lower bound and r.h.s.\ vectors of the refining
%% problem, whereas $c^{\mathrm{res}},l^{\mathrm{res}},b^{\mathrm{res}}$
%% represent their corresponding residuals for a given solution.

\paragraph{Cumulative dual LP.}
The corresponding dual problem is
\begin{equation}
\label{eq:DT}
\tag{$D_T$}
\left.
\begin{array}{rrll}
\max
& y_1+y_2 \\
\textnormal{s.t.:}
& y_1-3\displaystyle\sum_{\substack{j=1\\j\textnormal{ odd}}}^{T+1}y_j &\leq& 0, \\
& y_1-3\displaystyle\sum_{\substack{j=1\\j\textnormal{ even}}}^{T+1}y_j &\leq& 0, \\
& \epsilon^{-j}(y_2+y_{j+2}) &\leq& \epsilon^T+\epsilon^{-j} \qquad \forall j\in\{1,\ldots,T\}, \\
& \epsilon^{-T}(y_2+y_{T+3}) &\leq& \epsilon^T+\epsilon^{-T}, \\
& y_{j+2} &\leq& \zeroes \qquad \forall j\in\{1,\ldots,T+1\}, \\
& y_1,y_2 &\in& \mathbb{R}.
\end{array}
\right\}
\end{equation}

\paragraph{Elementary properties of \eqref{eq:PT}.}
It is easy to verify that the constraint matrix $A$ of \eqref{eq:PT} has full row rank $T+3$: the slack variables $s$ form an identity matrix in rows $3,\dots,T+3$, variables $w$ appear only in the first row, and the second row contains $z$ variables only; thus, it is straightforward to pick a basis. Furthermore, the problem admits an optimum, and therefore so does its dual.
\begin{proposition}[Optimal solutions of \eqref{eq:PT}]
  \label{prop:small-optimal-face}
  The cumulative primal LP \eqref{eq:PT} has optimal value $1+\epsilon^{2T}$. Every optimal solution $(w^\star, z^\star, s^\star)$ satisfies $z^\star_j=0 \; \forall j\in\{1,\ldots,T-1\}$, $z^\star_{T}=\lambda\epsilon^T, z^\star_{T+1}=(1-\lambda)\epsilon^T$ for some $\lambda\in[0,1]$, and $\|w^\star\|_\infty \le 1$, $\|z^\star\|_\infty \le 1$, $\|s^\star\|_\infty \le 3$.
\end{proposition}

\begin{proof}
Take any feasible solution to \eqref{eq:PT}. The second constraint reads
\begin{equation*}
  \sum_{j=1}^{T}\epsilon^{-j}z_j + \epsilon^{-T}z_{T+1} = 1.
\end{equation*}
Therefore the objective function satisfies
\begin{equation*}
  \sum_{j=1}^{T}
  \bigl(\epsilon^T+\epsilon^{-j}\bigr)z_j +
  \bigl(\epsilon^T+\epsilon^{-T}\bigr)z_{T+1} =
  1+\epsilon^T\sum_{j=1}^{T+1}z_j.
\end{equation*}
Moreover,
\begin{equation*}
  1 = \sum_{j=1}^{T}\epsilon^{-j}z_j + \epsilon^{-T}z_{T+1} \leq \epsilon^{-T}\sum_{j=1}^{T+1}z_j.
\end{equation*}
Hence $\sum_{j=1}^{T+1}z_j\geq\epsilon^T$, and consequently the objective function value is $\geq 1+\epsilon^{2T}$. It is straightforward to verify that the vector $z^\star$ as defined in the proposition statement satisfies the second constraint and attains objective function value $1+\epsilon^{2T}$. Furthermore, we need $\sum_{j=1}^{T+1}z_j = \epsilon^T$ to attain this value, and because only $z_T, z_{T+1}$ have coefficient $\epsilon^{-T}$ in the second constraint, that constraint implies $z_j=0 \; \forall j\in\{1,\ldots,T-1\}$, $z_T + z_{T+1} = \epsilon^T$ in every optimal solution. To complete the solution vector, one can choose $w^\star_1=\lambda, w^\star_2=1-\lambda$ if $T$ is odd, and $w^\star_1=1-\lambda, w^\star_2=\lambda$ if $T$ is even; there are other possible values for the $w$ variables. $s^\star$ is implied by $w^\star$ and $z^\star$.

The bounds on the infinity norms follow from the fact that all variables are nonnegative, and the first two constraints imply that no $w$ or $x$ variable can be larger than 1; the definition of the slack variables then implies that they cannot be larger than 3.
\end{proof}

%% At a high level, our argument has two layers. First, we study solvers that return primal basic, complementary-slack solutions and, at the degenerate vertices identified in Section \ref{sec:vertexstruct}, use a canonical basis completion. The canonical completion is a restriction on the basis reported by the solver at degenerate boundary singleton vertices. It does not change the primal vertex, but it does choose the associated basis dual and reduced cost vector. Thus it is not post-processing of a fixed primal-dual output: different degenerate bases can yield different dual infeasibilities and different classifications as lucky or non-lucky. Under this convention, if the returned solution is not unexpectedly more accurate than requested—what we later call lucky convergence—then the first active structural scale is forced to occur at the refinement index. This is the central obstruction in the paper: a non-lucky basic refining path cannot skip past the parity-switching scale. Second, we complement this structural result with a concrete algorithmic example, showing that a natural warm-started simplex implementation with Bland's rule realizes the corresponding exponentially-scaling path.

We can now state the main theorem of the paper. It shows that, for our LP family, an LP solver that returns basic, complementary-slack solutions is either \emph{lucky}, or it follows a path with exponential growth of the refining solutions. The concept of lucky convergence is formally defined in Definition~\ref{def:lucky-convergence}; essentially, a solution to the refinement problem is lucky if its dual infeasibility is $\bigO{\epsilon^2}$, which is better than the infeasibility that one would be able to guarantee with a limited-precision LP solver executed with tolerance $\epsilon$. Thus, a non-lucky solver returns a solution with dual infeasibility $\bigTh{\epsilon}$. The theorem statement also relies on the concept of a canonical basis. This is formally defined in Definition~\ref{def:canonicalbasis}, and it is our way of dealing with degenerate bases: a few vertices of the polytope (those with $w_1$ or $w_2$ equal to $0$) correspond to multiple degenerate bases, and we need a mechanism to pick one of these bases, because the dual infeasibility depends on it. We note that the assumption that the LP solver returns a canonical basis does not appear to be particularly strong: as we discuss in Section~\ref{sec:vertexstruct}, as $T$ gets large the fraction of nondegenerate vertices of the polytope goes to 1, and the unique basis for those vertices is canonical by definition.
\begin{theorem}[Exponential growth of refining solutions]
\label{thm:main-hardness}
Fix an integer $T\ge 3$ and a constant $0<\beta_T \leq 1$ that may depend polynomially on $T$. Define
\begin{equation*}
  \epsilon_0(T,\beta_T):=
  \min\left\{
  \frac{1}{2(T+1)},
  \frac{\beta_T}{T+1},
  \frac{1}{2}
  \right\}.
\end{equation*}
Let $\epsilon\in\mathbb{Q}\cap(0,\epsilon_0(T,\beta_T))$ and $0< \epsfin < \beta_T\epsilon^{T-2}$. Suppose that we apply Algorithm~\ref{alg:ir} to problem \eqref{eq:PT} with parameters $T, \epsilon, \epsfin$. Suppose further that the LP solver always returns a primal-feasible basic pair with a canonical basis, and that it does not get lucky, i.e., $\hat{\delta}_D^{(t)} \ge \beta_T\epsilon$ for all $t=1,\dots,T-2$, and $\delta_D^{(0)} \ge \beta_T\epsilon$ for the initial solve. Then, for every $t \in \{1,\ldots,T-2\}$ we have $\|\hat x^{(t)}\|_\infty = \bigOm{\epsilon^{-t}}$.
%% \begin{equation*}
%% \|\hat w^{(t)}\|_\infty \ge \frac{1}{3(T+1)}\epsilon^{-t}.
%% \end{equation*}
\end{theorem}
Our second main result shows that the exponential growth identified in Theorem~\ref{thm:main-hardness} occurs when we apply Algorithm~\ref{alg:ir} using a concrete implementation of the simplex algorithm, that relies on Bland's pivoting rule and warm-starts each LP solve in a natural way.
\begin{theorem}[Simplex realization of a sequence of solutions with exponentially-growing size]
\label{thm:main-simplex}
Fix an integer $T\ge3$, let $\epsilon\in\mathbb{Q}\cap(0,\frac12]$,
  and let $0<\epsfin<\frac34\epsilon^{T-2}$. Suppose we apply
  Algorithm~\ref{alg:ir} to problem \eqref{eq:PT} with parameters $T,
  \epsilon, \epsfin$, with an LP solver that implements the primal
  simplex algorithm in exact arithmetic, using Bland's pivoting rule,
  optimality tolerance $\epsilon$ (i.e., stopping when all reduced
  costs are greater than $-\epsilon$), and with the decision variables
  reordered as indicated in \eqref{eq:bland-order}. Suppose further
  that the solution returned at the initial solve
  (line~\ref{line:cumulprob} of Algorithm~\ref{alg:ir}) is the basic
  feasible solution corresponding to the basis
  \eqref{eq:initial-simplex-basis}, and every subsequent LP solve is
  warm-started from the final basis returned returned by the preceding
  solve. Then, for every $t\in\{1,\ldots,T-2\}$ we have $ \|\hat
  x^{(t)}\|_\infty =\bigOm{\epsilon^{-t}}$.
\end{theorem}
In the proof of Theorem~\ref{thm:main-simplex} we show that the
initial basis \eqref{eq:initial-simplex-basis} is $\epsilon$-optimal
in the stated manner. Thus, it is sufficient for the simplex
implementation to explore that basis at any time during the initial
solve, and it will terminate reporting it as sufficiently optimal.

%% Theorems~\ref{thm:main-hardness} and \ref{thm:simplexmain} serve different purposes. Theorem~\ref{thm:main-hardness} identifies the basis-resolved structural obstruction for basic solutions in the absence of lucky convergence, while Theorem~\ref{thm:simplexmain} shows that a natural and widely used LP method realizes the corresponding exponentially-growing path.

\section{Vertex structure of the polytope \eqref{eq:PT}}
\label{sec:vertexstruct}
To prove Theorems~\ref{thm:main-hardness} and \ref{thm:main-simplex}
it is convenient to characterize the vertices of
\eqref{eq:PT}. Introduce the scaled variables
\begin{equation*}
z'_j:=\epsilon^{-j}z_j,
\qquad
j=1,\ldots,T,
\qquad \qquad \text{and} \qquad \qquad
z'_{T+1}:=\epsilon^{-T}z_{T+1}.
\end{equation*}
Define
\begin{equation*}
O:=\{j\in\{1,\ldots,T+1\}:j\textnormal{ odd}\},
\;\;
E:=\{j\in\{1,\ldots,T+1\}:j\textnormal{ even}\},
\;\;
p(j) := \begin{cases} 1 & j \in O \\ 2 & j \in E.\end{cases}
\end{equation*}

In terms of the new variables, the feasible region of \eqref{eq:PT} can be written as
\begin{equation}
  \label{eq:PTp}
  \tag{$P'_T$}
P'_T := \left\{
(w_1,w_2,z'): w_1,w_2\geq 0, \; w_1+w_2=1, \;
\sum_{j=1}^{T+1}z'_j=1, \; 0\leq z'_j\leq 3w_{p(j)}
\right\}.
\end{equation}
Because the difference between $P'_T$ and the feasible region of \eqref{eq:PT} is only scaling of the $z$ variables, any result on the vertex structure of $P'_T$ immediately translates into a result on \eqref{eq:PT}. We initially study $P'_T$ because it simplifies the exposition, and later revert to \eqref{eq:PT}.

%We call \emph{active} a variable that is nonzero in a given solution.
Recall that $w$ are called controller variables and $z$ (or their scaled version $z'$) are called structural variables. It is immediate to observe that $P'_T \cap \{w_1 = 1\}$ is a simplex over the odd structural variables, and $P'_T \cap \{w_2 = 1\}$ is a simplex over the even structural variables. There are also vertices where both $w_1, w_2 > 0$; we characterize them next.
\begin{lemma}[Characterization of $z'$ at vertices of $P'_T$]
  \label{lem:vertexbinding}
  If $(w_1,w_2,z')\in P$ is a vertex of $P'_T$ with $w_1,w_2>0$, then $z'_j\in\{0,3w_{p(j)}\}$ for all $j=1,\ldots,T+1$.
\end{lemma}
\begin{proof}
  It is immediate to observe that $P'_T$ is $(T+1)$-dimensional (it has $T+3$ variables and two equality constraints). Because $w_1, w_2 > 0$, the lower bounding constraints for $w$ are not binding. Then we need $T+1$ binding constraints to define a vertex, so for each $z'_j$ either the lower or upper bound of $0\leq z'_j\leq 3w_{p(j)}$ must be binding.
\end{proof}

\begin{definition}[Support of a vector]
  \label{def:support}
  Let $z\in\mathbb{R}^{T+1}$, $z\geq \zeroes$. We define its \emph{support}, denoted
  $S(z)$, as
  \begin{equation*}
    S(z):= \{j\in\{1,\ldots,T+1\}:z_j>0\}.
  \end{equation*}
\end{definition}

\begin{proposition}[Vertices of $P'_T$ and \eqref{eq:PT}]
  \label{prop:vertexcases}
  Let $(w_1,w_2,z') \in P'_T$ be a vertex. Then one of the
  following holds:
  \begin{enumerate}
    \item[(1O).] $w_1=1$, $w_2=0$, $S(z') = \{j\}$ for some odd $j \in \{1,\dots,T+1\}$;
    \item[(1E).] $w_1=0$, $w_2=1$, $S(z') = \{j\}$ for some even $j \in \{1,\dots,T+1\}$;
    \item[(2O).] $w_1 = \frac{1}{3\ell}, w_2 = 1 - \frac{1}{3\ell}$ for some $\ell \in \{1,\dots,T+1\}$, $|S(z')| = \ell$ and $S(z')$ contains only odd indices;
    \item[(2E).] $w_1 = 1-\frac{1}{3k}, w_2 = \frac{1}{3k}$ for some $k \in \{1,\dots,T+1\}$, $|S(z')| = k$ and $S(z')$ contains only even indices.
  \end{enumerate}
\end{proposition}
\begin{proof}
  Cases 1O and 1E are trivial: $P'_T \cap \{w_1 = 1\}$ and $P'_T \cap \{w_2 = 1\}$ are simplices over odd-index and even-index $z'$ variables, so the vertices have a single nonzero $z'$ variable.

  To understand the structure of the vertices with $w_1, w_2 > 0$, note that by Lemma~\ref{lem:vertexbinding}, at any any such vertex there must exist  $S_E\subseteq E, S_O\subseteq O$ such that
  \begin{equation*}
    z'_j=
    \begin{cases}
      3w_1, & j\in S_O,\\
      3w_2, & j\in S_E,\\
      0, & j\notin S_E \cup S_O.
    \end{cases}
  \end{equation*}
  Let $\ell:=|S_O|, k:=|S_E|$. Then the equality constraints in \eqref{eq:PTp} imply that
  \begin{equation*}
    3\ell w_1+3k w_2=1, \qquad w_1+w_2=1,
  \end{equation*}
  and so $\ell \neq k$. Solving this system gives
  \begin{equation*}
    w_1=\frac{\frac{1}{3}-k}{\ell-k},
    \qquad
    w_2=\frac{\ell-\frac{1}{3}}{\ell-k}.
  \end{equation*}
  We are only interested in nonnegative solutions, and it is immediate to notice that if both $\ell > 0, k > 0$ the numerators of $w_1, w_2$ are of opposite sign, so they cannot be both nonnegative. This implies that one among $\ell, k$ must be zero. Case 2O corresponds to $\ell > 0, k = 0$ and case 2E corresponds to $\ell = 0, k > 0$.
\end{proof}

%% \begin{corollary}
%% Vertices of types $1$ and $2$ are boundary singleton vertices,
%% equivalently $|S(x)|=1$ and one controller variable is zero.
%% Vertices of types $3$ and $4$ with exactly one active structural
%% variable are interior singleton vertices. Vertices of types $3$ and
%% $4$ with at least two active structural variables are split
%% vertices.

%% At every interior vertex, the positive variables in the standard-form
%% representation determine the basis uniquely. Boundary singleton
%% vertices are degenerate and may admit multiple bases.
%% \end{corollary}

Proposition~\ref{prop:vertexcases} separates the vertices of $P'_T$, and therefore of \eqref{eq:PT}, into two types---we call them type 1 and type 2---and each type has odd and even vertices, abbreviated O or E. As our analysis is concerned with basic solutions, we now study the bases that correspond to these vertices. Vertices of type 2 correspond to a unique basis of \eqref{eq:PT}. On the contrary, at vertices of type 1 the standard-form representation is degenerate: one controller variable is zero, and several choices of zero-valued basic variables correspond to the same vertex. Because reduced costs are basis-dependent, and we use the reduced costs to measure dual infeasibility, we fix this ambiguity by picking a canonical basis for these degenerate cases. Lemmas~\ref{lem:interior-unique} and \ref{lem:boundary-expose-terminal} below show that canonical bases are indeed bases for \eqref{eq:PT}.
\begin{definition}[Canonical basis]
  \label{def:canonicalbasis}
  Let $x=(w,z,s)$ be a vertex of \eqref{eq:PT}. If $x$ is of type 1, let $S(z) = \{r\}$. Let $h\in\{T,T+1\}$, $h$ odd if $r$ is even, $h$ even if $r$ is odd. We define the canonical basis of $x$ as
  \begin{equation*}
    B(x):= \{w_1,w_2,z_r\} \cup \{s_j:j\in\{1,\ldots,T+1\} \setminus\{h\}\}.
  \end{equation*}
  If $x$ is of type 2, we define its canonical basis as
  \begin{equation*}
    B(x) := \{w_1,w_2\} \cup \{z_j:j\in S(z)\} \cup \{s_j:j\in\{1,\ldots,T+1\} \setminus S(z)\}.
  \end{equation*}
\end{definition}
\begin{lemma}[Type 2 vertices are nondegenerate]
  \label{lem:interior-unique}
  Let $x=(w,z,s)$ be a vertex of \eqref{eq:PT} of type 2. The unique basis of \eqref{eq:PT} that corresponds to $x$ is $B(x)$.
\end{lemma}
\begin{proof}
  By Proposition~\ref{prop:vertexcases}, all nonzero $z$ variables have the same parity and $w > 0$. For each $j\in S(z)$, $z_j$ is positive and $s_j$ is zero, because the associated upper bounding constraint is tight. For each $j\notin S(z)$, $z_j$ is zero and $s_j$ is positive, because $w_{p(j)}$ is positive. Therefore the positive variables at $x$ are precisely
  \begin{equation*}
    w_1,\ w_2,\ z_j\ (j\in S(z)),\ s_j\ (j\notin S(z)).
  \end{equation*}
  Their number is
  \begin{equation*}
    2+|S(z)|+\bigl((T+1)-|S(z)|\bigr)=T+3,
  \end{equation*}
  which is exactly the number of equality constraints in the standard-form problem \eqref{eq:PT}. So the corresponding basis is unique.
\end{proof}
\begin{lemma}[Type 1 vertices admit a canonical basis]
  \label{lem:boundary-expose-terminal}
  Let $x=(w,z,s)$ be a vertex of \eqref{eq:PT} of type 1. Let $S(z) =
  \{r\}$, and let $h\in\{T,T+1\}$, $h$ odd if $r$ is even, $h$ even if
  $r$ is odd. Then $B(x)$ is a basis representing $x$. Moreover, the
  nonbasic variable $z_h$ defines an edge that connects $x$ to the
  vertex $v = (w^\circ,z^\circ,s^\circ)$ of type 2 with
  $S(z^\circ) = \{h\}$.
\end{lemma}
\begin{proof}
  First, note that $|B(x)| = T+3$, which is the number of constraints
  of \eqref{eq:PT}. We prove that the submatrix $A_{B(x)}$ of \eqref{eq:PT} is
  nonsingular, by looking at solutions of the system $A_{B(x)} \lambda =
  \zeroes$. 

  In the second row of \eqref{eq:PT}, $\sum_{j=1}^{T}\epsilon^{-j}z_j
  +\epsilon^{-T}z_{T+1}$, only $z_r$ is basic, so $\lambda_{z_r} =
  0$. In the row $-3w_{p(h)}+\epsilon^{-j}z_h+s_h$ only $w_{p(h)}$ is
  basic, so $\lambda_{w_{p(h)}} = 0$. Then the first row $w_1 + w_2$
  implies that $\lambda_{w_1} = \lambda_{w_2} = 0$. Finally, in all
  other rows only the slack $s_j, j \in \{1,\dots,T+1\} \setminus
  \{h\}$, is basic, implying that $\lambda_{s_j} = 0$. Thus, the only
  solution to $A_{B(x)} \lambda = \zeroes$ is $\lambda = \zeroes$,
  proving that $B(x)$ is a basis. It represents $x$ because given
  $B(x)$, only the corresponding type 1 vertex has nonzero support
  over $B(x)$.

  We now look at the edge defined by the nonbasic variable
  $z_h$. Suppose, without loss of generality, that $r$ is even and $h$
  is odd. Recall that $h \in \{T,T+1\}$. Because $s_h$ is nonbasic and
  therefore equal to $0$, the row containing $s_h$ reads
  $\epsilon^{-T}z_h=3w_{p(h)}$. At $z_h =\epsilon^T$ this equation
  implies $w_{p(h)} = 1/3$, and the second row of \eqref{eq:PT}, which
  is the only other row containing $z_h$, implies $z_r = 0$. In terms
  of pivots, $z_h$ enters the basis, $z_r$ leaves. These values and
  basis define another vertex of type 2, with opposite parity compared
  to $x$, so this has to be the endpoint of the edge.  
\end{proof}
For crucial parts of the paper we make the assumption that the LP
solver used in Algorithm~\ref{alg:ir} returns a canonical basis. We
argue that this is not a particularly strong assumption, for two
reasons. The first reason is that as $T$ gets large, the overwhelming
majority of vertices are nondegenerate, and for nondegenerate vertices
the canonical basis assumption is automatically satisfied. Indeed, using
Proposition~\ref{prop:vertexcases}, we observe that there are
$\bigO{T}$ vertices of type 1 (due to the two
$\frac{T+1}{2}$-dimensional simplices), but $\bigOm{2^{T}}$ vertices
of type 2 (for each possible value of $|S(z')|$, we sum over all
subsets of $O, E$ of that size; and $\sum_{k=1}^{n} \binom{n}{k} =
2^n-1$) for which the basis is unique and canonical. The second reason
is that the purpose of this paper is to provide evidence that it is
difficult to provably avoid the exponential growth of the refining
solutions; to get around our result, one would have to show that the
LP solver sufficiently often does not return a canonical basis, and
also that exponential growth cannot occur with other types of
bases. As remarked in Section~\ref{sec:intro}, our paper provides a
blueprint for showing the exponential lower bound that may not hold in
all cases, but that can plausibly be adapted to specific LP solvers.

We conclude this section by characterizing the reduced costs
associated with canonical bases for vertices of type 1. Reduced costs
are crucial to our analysis, as they determine the dual infeasibility.
\begin{lemma}[Characterization of reduced costs for type 1 vertices]
  \label{lem:boundary-reduced-costs}
  Let $x = (w,z,s)$ be a vertex of \eqref{eq:PT} of type 1. Let $S(z) =
  \{r\}$. Then for every nonbasic structural variable $z_j, j \neq r$,
  its reduced cost with respect to the basis $B(x)$ is
  \begin{equation*}
    \bar c^{B(x)}_{z_j} = \epsilon^T \left(1-\epsilon^{\min\{r,T\}-\min\{j,T\}}\right).
  \end{equation*}
\end{lemma}
\begin{proof}
  Let $h\in\{T,T+1\}$, $h$ odd if $r$ is even, $h$ even if $r$ is
  odd. The dual vector associated with $B(x)$ is defined by the linear
  system $A_{B(x)}^\top\pi=c_{B(x)}$, where $c_{B(x)}$ are the
  components of the objective vector corresponding to $B(x)$.  We
  assume that the order of the variables in $B(x)$, and therefore
  $A_{B(x)}$, is exactly as written in Definition~\ref{def:canonicalbasis}, i.e.,
  $\{w_1,w_2,z_r\} \cup \{s_j:j\in\{1,\ldots,T+1\}
  \setminus\{h\}\}$. Note that $\pi \in \mathbb{R}^{|T|+3}$, and
  because of the duality relationship its elements correspond to rows
  of $A$.

  For every $j \in \{1,\dots,T+1\} \setminus \{h\}$, $s_j \in B(x)$,
  $A_{s_j}$ is an identity column and $c_{s_j} = 0$, so $\pi_{j+2} =
  0$. Also, $w_1, w_2 \in B(x)$ and $c_{w_1} = c_{w_2} = 0$. Because
  $h$ has parity opposite to $r$, together with $\pi_{j+2} = 0$ for $j
  \in \{1,\dots,T+1\} \setminus \{h\}$, the equation of
  $A_{B(x)}^\top\pi=c_{B(x)}$ corresponding to $w_{p(r)}$ gives
  $\pi_{1}=0$. Then the equation corresponding to the other $w$
  variable implies $\pi_{h+2}=0$. Thus, $\pi_j = 0$ for $j \in \{1\}
  \cup \{3,\dots,T+3\}$.

  Now we look at the equation for $z_r$. Column $A_{z_r}$ has
  coefficient $\epsilon^{-\min\{r,T\}}$ in the second row and in the row that
  contains $s_r$ (the row that upper bounds $z_r$). Because $\pi_{r+2}
  = 0$, this equation gives $\epsilon^{-\min\{r,T\}} \pi_2 =
  \epsilon^T+\epsilon^{-\min\{r,T\}}$, which implies $\pi_2 = 1+\epsilon^{T+\min\{r,T\}}$.
  Then for every nonbasic $z_j$ we have
  \begin{equation*}
    \begin{aligned}
      \bar c^{B(x)}_{z_j}
      &= \epsilon^T+\epsilon^{-\min\{j,T\}} - \epsilon^{-j}\pi_2 
      = \epsilon^T+\epsilon^{-\min\{j,T\}} - \epsilon^{-\min\{j,T\}}
      \left(1+\epsilon^{T+\min\{r,T\}} \right) \\
      &= \epsilon^T-\epsilon^{T+\min\{r,T\}-\min\{j,T\}} 
      = \epsilon^T \left(1-\epsilon^{\min\{r,T\}-\min\{j,T\}}\right). 
    \end{aligned}
    \qedhere
  \end{equation*}
\end{proof}

\section{Proof of Theorem~\ref{thm:main-hardness}}
\label{sec:lucky-convergence}

The section proceeds as follows. We first define what we mean by
\emph{lucky convergence}, in terms of the dual infeasibility of the
cumulative solutions $(x^{(t)}, y^{(t)})$ constructed by
Algorithm~\ref{alg:ir}. Next, we relate reduced costs of the canonical
bases (Definition~\ref{def:canonicalbasis}) to normalized objective drops
along nonbasic directions. Finally, using results from
Section~\ref{sec:vertexstruct}, we show that unless the LP solver is
lucky, the index of the first nonzero $z$ variable in a cumulative
solution advances by one in every IR iteration, and this causes the
$w$ variables to bounce. The relationship between cumulative and
refining solutions of
Proposition~\ref{prop:cumulative-refining-correspondence} then
concludes the proof.

\subsection{Lucky convergence}
We use the largest dual infeasibility $\delta^{(t)}_D$ to measure progress of the IR algorithm (the sequence of solutions that we study is always primal feasible). The promise of Algorithm~\ref{alg:ir} is that $\delta_D^{(t)}=\bigO{\epsilon^{t+1}}$, see Theorem~\ref{thm:irconvergence}. We call the iteration lucky when $\delta_D^{(t)}$ decreases faster than promised by Theorem~\ref{thm:irconvergence}.

\begin{definition}[Lucky and non-lucky iterations]
  \label{def:lucky-scales}
  Let $T$ be fixed. Choose constants $\mu_T, \beta_T$, depending only on $T$. For every $t=1,\dots,T-1$, we say that iteration $t$ of Algorithm~\ref{alg:ir}, and the cumulative iterate $(x^{(t)}, y^{(t)})$ after iteration $t$, are:
  \begin{enumerate}
  \item \emph{on the intended scale} if $\delta_D^{(t)} \le \mu_T\epsilon^{t+1}$;
  \item \emph{lucky} if $\delta_D^{(t)} \le \mu_T\epsilon^{t+2}$;
  \item \emph{non-lucky} if $\delta_D^{(t)} \ge \beta_T\epsilon^{t+1}$.
  \end{enumerate}
\end{definition}
To be precise, one should indicate the constants $\mu_T, \beta_T$ when
using the concepts of Definition~\ref{def:lucky-scales}. However,
because in this paper we are concerned with asymptotic statements that
hold when $\epsilon$ is sufficiently small, we can think of $\mu_T,
\beta_T$ simply as some constants that do not depend on $\epsilon$:
their exact value is immaterial. Thus, a cumulative primal-dual pair
$(x^{(t)}, y^{(t)})$ of Algorithm~\ref{alg:ir} is on the intended
scale, but non-lucky, if $\delta^{(t)}_D = \bigTh{\epsilon^{t+1}}$.
\begin{definition}[Lucky convergence]
  \label{def:lucky-convergence}
  We say that the LP solver exhibits \emph{lucky convergence} at iteration $t$ if iteration $t$ of Algorithm~\ref{alg:ir} is lucky in the sense of Definition~\ref{def:lucky-scales}.
\end{definition}
Equivalently, Definition~\ref{def:lucky-scales} can be stated in terms
of the dual infeasibility of the refining problem
$\widehat{\delta}^{(t)}_D$. Proposition~\ref{prop:cumulative-refining-correspondence}
immediately implies that an iteration is on the intended scale if
$\widehat{\delta}^{(t)}_D = \bigO{\epsilon}$, lucky if
$\widehat{\delta}^{(t)}_D = \bigO{\epsilon^2}$, and non-lucky if
$\widehat{\delta}^{(t)}_D = \bigOm{\epsilon}$. Thus, an iteration of
Algorithm~\ref{alg:ir} is lucky if the LP solver returns a solution to
the refining problem that is significantly better than the optimality
tolerance $\epsilon$ with which it is executed. We call this ``lucky''
because, although it can happen in practice, it is an event that we
cannot rely on.

\subsection{The normalized objective and its properties}
To simplify certain calculations with reduced costs, we work with a normalized version of the objective function of \eqref{eq:PT}.
\begin{definition}[Normalized objective]
  \label{def:normalizedobj}
  Let $x = (w,z,s)$ be feasible for \eqref{eq:PT}. We define its \emph{normalized objective} $f(v)$ as
  \begin{equation*}
    f(v) := \epsilon^{-T}(c^\top v-1) = \epsilon^{-T}
    \left(\sum_{j=1}^{T} \bigl(\epsilon^T+\epsilon^{-j}\bigr)z_j +
    \bigl(\epsilon^T+\epsilon^{-T}\bigr)z_{T+1} - 1 \right).
  \end{equation*}
\end{definition}
\begin{lemma}[Expression for the normalized objective]
  \label{lem:normalizedobjsum}
  Let $x = (w,z,s)$ be feasible for \eqref{eq:PT}. Then $f(v) = \sum_{j=1}^{T+1}z_j$.
\end{lemma}
\begin{proof}
  Because $x$ is feasible, the second constraint of \eqref{eq:PT}
  gives $\sum_{j=1}^{T}\epsilon^{-j}z_j + \epsilon^{-T}z_{T+1} =
  1$. So
  \begin{equation*}
    c^\top v = \sum_{j=1}^{T} \bigl(\epsilon^T+\epsilon^{-j}\bigr)z_j +
    \bigl(\epsilon^T+\epsilon^{-T}\bigr)z_{T+1} = 1 + \epsilon^T\sum_{j=1}^{T+1}z_j. \qedhere
  \end{equation*}
\end{proof}
We also work with the nonbasic directions defined by a basis. Their definition is standard in the literature; we recall it below for the sake of self-containment.
\begin{definition}[Nonbasic direction]
  \label{def:nonbasicdir}
  Let $(B, N)$ be a partition of the decision variables into basic and nonbasic. Let $x_j \in N$. The nonbasic direction corresponding to $x_j$, denoted $d^{B,x_j}$, is defined as
  \begin{equation*}
    d^{B,x_j}_{x_j} = 1,\qquad d^{B,x_j}_N = \zeroes, \qquad d^{B,x_j}_B=-A_B^{-1}A_{x_j}.
  \end{equation*}
\end{definition}
\begin{lemma}[Reduced costs as normalized edge slopes]
  \label{lem:reduced-costs-edge-slopes}
  Let $B$ be a basis for \eqref{eq:PT}, and let $x$ be the corresponding
  vertex. Let $x_j \not\in B$ and suppose that
  \begin{equation*}
    v = x + \theta d^{B,x_j} \qquad \text{for some } \theta > 0,
  \end{equation*}
  for some other vertex $v$, i.e., the nonbasic direction
  corresponding to $x_j$ leads to a distinct vertex with positive step
  length $\theta$. Then
  \begin{equation*}
    -\bar c^B_{x_j} = \epsilon^T\frac{f(x)-f(v)}{\theta}.
  \end{equation*}
\end{lemma}
\begin{proof}
  It is well known that the rate of change of the objective function
  of an LP along a nonbasic direction is equal to the reduced cost of
  the corresponding nonbasic variable, and this relationship holds at
  equality until the basis changes. Then, using
  Definition~\ref{def:normalizedobj}, it immediately follows that 
  \begin{equation*}
    \frac{d}{d\theta}f(x+\theta d^{B,x_j}) = \epsilon^{-T}\bar c^B_{x_j}.
  \end{equation*}
  By linearity of $f$, 
  \begin{equation*}
    f(v) = f(x) + \theta \frac{d}{d\theta}f(x+\theta d^{B,x_j}) = f(x) + \theta \epsilon^{-T}\bar c^B_{x_j}.
  \end{equation*}
  Rearranging the terms proves the statement.
\end{proof}
Lemma~\ref{lem:reduced-costs-edge-slopes} shows that we can compute
the reduced costs by using the normalized objective function $f$
divided by edge lengths. By Lemma~\ref{lem:normalizedobjsum}, $f$
admits a simple expression. This is our approach to characterize the
reduced costs and therefore the dual infeasibility of the vertices of
\eqref{eq:PT}.

As it turns out, the normalized objective of vertices is essentially
characterized by the index of their first nonzero structural
variable. Because this quantity is used frequently in the remainder of
the paper, we assign it its own symbol.
\begin{definition}[First nonzero variable]
  \label{def:qv}
  Let $x = (w,z,s)$ be a vertex of \eqref{eq:PT}. We denote the index of its first nonzero $z$ variable by $q(x)$, defined as
  \begin{equation*}
    q(x) := \min \left\{ \min\{ j \in S(z)\}, T\right\}.
  \end{equation*}
\end{definition}
In Definition~\ref{def:qv} we clip the value of $q(x)$ to the interval $[1,
  T]$ because the decision variable $z_{T+1}$ has the same
coefficients (in the objective and constraints) as $z_{T}$, i.e.,
$\epsilon^{-T}$. Indeed, $z_{T+1}$ is essentially a copy of $z_T$ that
we use to deal with odd/even asymmetries, and its effect on the
objective function value is the same, hence the definition of
$q(x)$. The next proposition makes precise the way in which $q(x)$
characterizes the normalized objective of $x$.

\begin{proposition}[Relationship between normalized objective and $q(x)$]
  \label{prop:normalized-objective}
  Let $x=(w,z,s)$ be a vertex of \eqref{eq:PT}. Then $f(x) =
  \frac{1}{|S(z)|} \sum_{j \in S(z)} \epsilon^{\min\{j,T\}} =
  \bigTh{\epsilon^{q(x)}}$, where the hidden constants may depend
  polynomially on $T$.
\end{proposition}

\begin{proof}
  By Proposition~\ref{prop:vertexcases}, $S(z)$ contains only odd or
  even indices. By Lemma~\ref{lem:normalizedobjsum}, $f(x) =
  \sum_{j=1}^{T+1}z_j$.

  For vertices of type 1, $|S(z)| = 1$ and the second constraint of
  \eqref{eq:PT} implies that $f(x) = \frac{1}{|S(z)|} \allowbreak
  \sum_{j \in S(z)} \epsilon^{\min\{j,T\}}$. For vertices of type 2,
  the upper bounding constraints that are binding at $x$ imply that,
  if we use the rescaled variables $z'$ of
  Proposition~\ref{prop:vertexcases}, $z'_j = \frac{1}{|S(z)|}$ for
  all $j \in S(z)$. Thus, translating back to the unscaled variables,
  $f(x) = \frac{1}{|S(z)|} \sum_{j \in S(z)}
  \epsilon^{\min\{j,T\}}$. So this expression holds for all vertices.

  By Definition~\ref{def:qv}, we have
  \begin{equation*}
    f(x) = \frac{1}{|S(z)|} \sum_{j \in S(z)} \epsilon^{\min\{j,T\}} \geq \frac{1}{|S(z)|}\epsilon^{q(x)} \geq \frac{1}{T+1}\epsilon^{q(x)},
  \end{equation*}
  because $|S(z)| \le \ceil{(T+1)/2}$. On the other hand, $q(x) \le j$ for every $j\in S(z)$, so $\epsilon^{j} \le \epsilon^{q(x)}$. Therefore
  \begin{equation*}
    f(x) = \frac{1}{|S(z)|} \sum_{j \in S(z)} \epsilon^{\min\{j,T\}}
    \leq \epsilon^{q(x)}. \qedhere
  \end{equation*}  
\end{proof}
Let $x = (w,z,s)$ be a vertex of \eqref{eq:PT}. As a direct consequence of Proposition~\ref{prop:normalized-objective}, for a given refinement iteration $t$, we have the following implications:
\begin{align*}
  q(x)\leq t \quad\Longrightarrow&\quad f(x)=\bigOm{\epsilon^t}, \\
  q(x)=t+1 \quad\Longrightarrow&\quad f(x)=\bigTh{\epsilon^{t+1}}, \\
  q(x)\geq t+2 \quad\Longrightarrow&\quad f(x)=\bigO{\epsilon^{t+2}}.
\end{align*}
This points to our strategy for the proof of
Theorem~\ref{thm:main-hardness}: we relate $q(x)$ to the quality of an
IR iteration, as defined in Definition~\ref{def:lucky-scales}. If $q(x)$ is
too small compared to $t$ an iteration would not make enough progress,
and if it is too large the iteration is lucky. We analyze this aspect
in the next subsection.

\subsection{Relationship between $q(x)$ and luckiness}
The previous section relates $q(x)$ to the normalized objective. We
now extend the analysis to relate $q(x)$ to the level of success of an
IR iteration, i.e., whether vertex $x$ represents enough progress at
some iteration $t$ or not. To do so, we need to bound the length of
certain edges of the polytope \eqref{eq:PT}, so that we can use
Lemma~\ref{lem:reduced-costs-edge-slopes}.

\begin{lemma}[Bounds on the edge lengths at canonical bases]
  \label{lem:edge_length-bound}
  Let $x=(w,z,s)$ be a vertex of \eqref{eq:PT}. Let $B(x)$ be its
  canonical basis. Let $x_j \not\in B(x)$ and suppose that $v = x +
  \theta d^{B(x),x_j}$ is a distinct vertex of \eqref{eq:PT} for some
  $\theta > 0$. Then
  \begin{equation*}
    \theta\geq\frac{1}{T+1}\epsilon^T.
  \end{equation*}
  Furthermore, if $x_j \in \{z_{T},z_{T+1}\}$, then $\theta\leq\epsilon^T$.
\end{lemma}
\begin{proof}
  Let $v=(w^\circ,z^\circ,s^\circ)$. By Definition~\ref{def:nonbasicdir},
  $d^{B(x),x_j}_{x_j}=1$. Because $x_j \not\in B(x)$, $x_j = 0$. Thus,
  $\theta$ is precisely the value $v_j$.

  First, suppose $x_j$ is $z_h$ for some $h \in \{1,\dots,T+1\}$. By
  Proposition~\ref{prop:normalized-objective}, if $z^\circ_h > 0$ then
  $z^\circ_h= \frac{\epsilon^{\min\{h,T\}}}{|S(z^\circ)|}$. Using
  $|S(z^\circ)|\leq T+1$, we obtain
  \begin{equation*}
    \theta=z^\circ_h \geq \frac{1}{T+1}\epsilon^{\min\{h,T\}} \geq
    \frac{1}{T+1}\epsilon^T.
  \end{equation*}
  If $h \in \{T,T+1\}$, then 
  \begin{equation*}
    \theta=z^\circ_h = \frac{\epsilon^T}{|S(z^\circ)|} \leq
    \epsilon^T.
  \end{equation*}
  
  Next, consider the case in which $x_j$ is $s_h$ for some $h \in
  \{1,\dots,T+1\} \setminus S(z)$. At vertices of type 2, using
  Definition~\ref{def:canonicalbasis} and Proposition~\ref{prop:vertexcases}
  we see that every positive slack is either equal to
  $\frac{1}{|S(z)|}$ or to $3-\frac{1}{|S(z)|}$. At vertices of type
  1, every positive slack is either $2$ or $3$. Thus $s^\circ_h \ge
  \frac{1}{T+1}$, and therefore $\theta\geq\frac{1}{T+1}$.  Because
  $0<\epsilon<1$, this implies
  \begin{equation*}
    \theta\geq\frac{1}{T+1}\epsilon^T. \qedhere
  \end{equation*}
\end{proof}
We provide a sharper characterization when the nonbasic variable is
the one among $z_T$ or $z_{T+1}$ which is of parity opposite to that
of the current vertex.
\begin{lemma}[Characterization of the nonbasic edge leading to the optimal face]
    \label{lem:terminal-edge}
    Let $x=(w,z,s)$ be a vertex of \eqref{eq:PT} of type 1O or 2O. Let $h \in \{T,T+1\}$, $h$ even. Then $v = x + \epsilon^T d^{B(x),z_h}$ is a vertex with a single nonzero $z$ variable, and $f(v) = \epsilon^T$. The symmetric statement exchanging ``odd'' and ``even'' also holds.
\end{lemma}
\begin{proof}
  Let $v=(w^\circ,z^\circ,s^\circ)$. Note that $z_h \not\in B(x)$.
  
  If $x$ is of type 1O, the proof follows directly from
  Lemma~\ref{lem:boundary-expose-terminal}, and the observation that
  $z^\circ_h = \epsilon^T$.

  Now suppose $x$ is of type 2O. Recall that $h$ is even and in a vertex of \eqref{eq:PT}, either odd or even $z_j$ can be nonzero---not both. Because
  \begin{equation*}
    B(x) = \{w_1,w_2\} \cup \{z_j:j\in S(z)\} \cup \{s_j:j\in\{1,\ldots,T+1\} \setminus S(z)\},
  \end{equation*}
  and $S(z)$ contains only odd indices, if $z_h$ increases all
  currently nonzero $z_j$ must go to zero. Then $S(z^\circ) = \{h\}$,
  and therefore, $z^\circ_h = \epsilon^T$. (In terms of simplex pivot,
  $z_h$ enters and there is a tie for the leaving variable, as $w_1$
  and $\{z_j:j\in S(z)\}$ all go to zero simultaneously.)

  In both cases, $z^\circ_h = \epsilon^T$. This implies that we moved
  by a step length $\theta = \epsilon^T$ and that $f(v) =
  \epsilon^T$. The proof for the symmetric statement follows the same
  structure.
\end{proof}
We now have all the tools to show two results about $q(x)$. The first result states that if $q(x)$ is too small compared to the iteration index $t$, the dual infeasibility of $x$ is not on the intended scale. The second states that $q(x)$ is too small compared to $t$, iteration $t$ is lucky.
\begin{lemma}[Vertices with small $q(x)$ have large dual infeasibility]
  \label{lem:low-level}
  Let $T$ be an integer, and let
  $0<\epsilon\le\frac{1}{2(T+1)}$. Suppose Algorithm~\ref{alg:ir} is
  executed on problem \eqref{eq:PT} with the above values of $T$ and
  $\epsilon$. For $t\in\{1,\ldots,T-1\}$, let $x=(w,z,s)$ be the
  vertex of \eqref{eq:PT} corresponding to the primal-dual pair
  $(x^{(t)}, y^{(t)})$ at iteration $t$, and suppose the LP solver
  reports the canonical basis $B(x)$ as the corresponding basis. If
  $q(x)\le t$, then
  \begin{equation*}
    \delta_D^{(t)} \ge
    \frac{1}{2(T+1)}\epsilon^t = \bigOm{\epsilon^t}.
  \end{equation*}
\end{lemma}
\begin{proof}
  Suppose $x$ is of type 1O or 2O; the proof if it is 1E or 2E is
  symmetric. Let $h \in \{T,T+1\}$, $h$ even. By
  Lemma~\ref{lem:terminal-edge}, $v = x + \epsilon^T d^{B(x),z_h}$ is
  a vertex with $f(v)=\epsilon^T$.

  Because $q(x)\leq t$, Proposition~\ref{prop:normalized-objective}
  implies that $f(x) \geq \frac{1}{T+1}\epsilon^t$. Because $t\le
  T-1$, $f(v) \le \epsilon^{t+1}
  \le\frac{1}{2(T+1)}\epsilon^t$. Combining these bounds, we obtain
  \begin{equation*}
    f(x)-f(v)\ge\frac{1}{2(T+1)}\epsilon^t.
  \end{equation*}

  Using Lemma~\ref{lem:reduced-costs-edge-slopes} and
  $\theta=\epsilon^T$, we have
  \begin{equation*}
    -\bar c_h^{B(x)} = \epsilon^T\frac{f(x)-f(v)}{\theta}
    = f(x)-f(v) \geq \frac{1}{2(T+1)}\epsilon^t = \bigOm{\epsilon^t}.
  \end{equation*}
  By definition, $\delta_D^{(t)} \ge -\bar c_h^{B(x)}$. This concludes
  of the proof.
\end{proof}

\begin{lemma}[Vertices with large $q(x)$ have small dual infeasibility]
  \label{lem:high-level}
  Let $T$ be an integer, and let
  $0<\epsilon\le\frac{1}{2(T+1)}$. Suppose Algorithm~\ref{alg:ir} is
  executed on problem \eqref{eq:PT} with the above values of $T$ and
  $\epsilon$. For $t\in\{0,1,\ldots,T-1\}$, let $x=(w,z,s)$ be a 
  vertex of \eqref{eq:PT} corresponding to the primal-dual pair
  $(x^{(t)}, y^{(t)})$ at iteration $t$, and suppose the LP solver
  reports the canonical basis $B(x)$ as the corresponding basis. If
  $q(x)\geq t+2$, then
  \begin{equation*}
    \delta_D^{(t)} \le (T+1)\epsilon^{t+2} = \bigO{\epsilon^{t+2}}.
  \end{equation*}
\end{lemma}

\begin{proof}
  If $x$ is optimal for \eqref{eq:PT} the proof is trivial, so from
  now on assume $x$ is not optimal. By
  Proposition~\ref{prop:normalized-objective}, $f(x) \leq
  \epsilon^{q(x)} \leq \epsilon^{t+2}$.
  
  Suppose first that $x$ is a type 2 vertex. By
  Lemma~\ref{lem:interior-unique}, $B(x)$ is nondegenerate. Hence,
  every edge corresponding to a nonbasic variable leads to a distinct
  adjacent vertex. Let $x_j \not\in B(x)$ and $\bar c_{x_j}^{B(x)}<0$.  Let
  $v=x+\theta d^{B(x),x_j}$ be a vertex. Since the edge is improving,
  Lemma~\ref{lem:reduced-costs-edge-slopes} gives
  \begin{equation*}
    -\bar c_{x_j}^{B(x)}
    =
    \epsilon^T\frac{f(x)-f(v)}{\theta}.
  \end{equation*}
  Moreover, $f(v)\geq 0$, so $f(x)-f(v)\leq f(x)$.
  By Lemma~\ref{lem:edge_length-bound}, $\theta\geq\frac{1}{T+1}\epsilon^T$.
  Therefore
  \begin{equation*}
    -\bar c_{x_j}^{B(x)} \leq \epsilon^T \frac{f(x)}{\epsilon^T/(T+1)}
    = (T+1)f(x) \leq (T+1)\epsilon^{t+2}.
  \end{equation*}
  
  Now suppose that $x$ is type 1 vertex and $S(z) = \{r\}$. Note that
  $s_h \in B(x)$ and it is the only nonbasic slack. Using the
  characterization of $\pi$ in the proof of
  Lemma~\ref{lem:boundary-reduced-costs}, the reduced cost of $s_h$ is
  zero. Thus, any negative reduced cost is associated with some $z_j
  \not\in B(x)$. By Lemma~\ref{lem:boundary-reduced-costs}, $\bar
  c_{z_j}^{B(x)} = \epsilon^T
  \left(1-\epsilon^{\min\{r,T\}-\min\{j,T\}} \right).$
  %If $\bar c_{z_j}^{B(x)}<0$, $j > r$.
  Therefore
  \begin{equation*}
    -\bar c_{z_j}^{B(x)}
    = \epsilon^T \left( \epsilon^{\min\{r,T\}-\min\{j,T\}}-1 \right)
    \leq \epsilon^{T+\min\{r,T\}-\min\{j,T\}}
    \leq \epsilon^{\min\{r,T\}}.
  \end{equation*}
  Also, $q(x) = \min\{r,T\}$. Thus $q(x)\geq t+2$ implies
  $\min\{r,T\}\geq t+2$, and hence $-\bar c_{z_j}^{B(x)} \leq
  \epsilon^{t+2}$.
  
  Thus, for every $x_j \not\in B(x)$ with $\bar c_{x_j}^{B(x)}<0$, we have
  $-\bar c_{x_j}^{B(x)} \le (T+1)\epsilon^{t+2}$. Then
  \begin{equation*}
    \delta_D^{(t)} = \max_{x_j \not\in B(x) :\bar c_{x_j}^{B(x)}<0}
    \{-\bar c_{x_j}^{B(x)}\}
    \le (T+1)\epsilon^{t+2} = \bigO{\epsilon^{t+2}}. \qedhere
  \end{equation*}
\end{proof}

The next proposition combines Lemmas~\ref{lem:low-level}
and~\ref{lem:high-level} to identify the ``right'' $q(x)$ at non-lucky
iterations.
\begin{proposition}[Progress in $q(x)$ at non-lucky iterations]
  \label{prop:forced-level}
  Let $T$ be an integer, let $\beta_T > 0$ (possibly depending on
  $T$), and let $0<\epsilon< \min\left\{ \frac{1}{2(T+1)},
  \frac{\beta_T}{T+1} \right\}$. Suppose Algorithm~\ref{alg:ir} is
  executed on problem \eqref{eq:PT} with the above values of $T$ and
  $\epsilon$. For $t\in\{1,\ldots,T-1\}$, let $x=(w,z,s)$ be a vertex
  of \eqref{eq:PT} corresponding to the primal-dual pair $(x^{(t)},
  y^{(t)})$ at iteration $t$, and suppose the LP solver reports the
  canonical basis $B(x)$ as the corresponding basis. Suppose further
  that $\delta_D^{(t)}\le\epsilon^{t+1}$ and
  $\delta_D^{(t)}\ge\beta_T\epsilon^{t+1}$.  Then $q(x)=t+1$.
  %% Consequently, if $t\leq T-2$, then the structural variable
  %% $z_{t+3}$ is active at $x$. If $t=T-1$, then at least one
  %% terminal structural variable, $z_{T+2}$ or $z_{T+3}$, is active.
\end{proposition}

\begin{proof}
  If $q(x)\leq t$, Lemma~\ref{lem:low-level} shows that
  $\delta_D^{(t)} \geq \frac{1}{2(T+1)}\epsilon^t$.
  Because $\epsilon<\frac{1}{2(T+1)}$, this would
  contradict $\delta_D^{(t)}\leq\epsilon^{t+1}$. Thus, $q(x) \ge t+1$.

  If $q(x)\geq t+2$, Lemma~\ref{lem:high-level} shows that
  $\delta_D^{(t)} \leq (T+1)\epsilon^{t+2}$.  Because
  $\epsilon<\frac{\beta_T}{T+1}$, this would contradict
  $\delta_D^{(t)}\geq\beta_T\epsilon^{t+1}$. Thus, $q(x) \le t+1$.
\end{proof}

\subsection{Putting it together}
We now complete the proof of Theorem~\ref{thm:main-hardness}. Note that because the dual infeasibility is zero at an optimal basic primal-dual pair, the final iteration of Algorithm~\ref{alg:ir} is necessarily lucky according to Definition~\ref{def:lucky-scales}, so we focus on indices $t \in \{1,\dots,T-2\}$.

\begin{proof}[Proof of Theorem~\ref{thm:main-hardness}]
  Using our convention for the names of the variables, we let $x^{(t)}
  = (w^{(t)},z^{(t)},s^{(t)})$ and $\hat{x}^{(t)} =
  (\hat{w}^{(t)},\hat{z}^{(t)},\hat{s}^{(t)})$.
  
  By assumption, every refining solution $(\hat{x}^{(t)},
  \hat{y}^{(t)})$ is a primal-feasible basic pair for the refining
  problem. By
  Proposition~\ref{prop:cumulative-refining-correspondence}, every
  $(x^{(t)}, y^{(t)})$ is also a primal-feasible basic pair for
  \eqref{eq:PT}. By Definition~\ref{def:lucky-scales} and the fact that the
  LP solver is not lucky, $\delta_D^{(t)} \geq \beta_T \epsilon^{t}$.
  By Theorem~\ref{thm:irconvergence}, $\delta_D^{(t)} \leq
  \epsilon^{t+1}$.

  We first show that if the LP solver does not get lucky, we execute
  at least $T-1$ refinement iterations. By assumption, primal
  feasibility and complementary slackness always hold. Hence, only the
  dual infeasibility condition can stop the algorithm.  For every
  $\tau \in\{1,\ldots,T-2\}$, the fact that the LP solver is not lucky
  at iteration $\tau-1$ implies
  \begin{equation*}
    \delta_D^{(\tau-1)} \geq \beta_T\epsilon^\tau \geq \beta_T\epsilon^{T-2}
    > \epsilon_{\mathrm{final}}.
  \end{equation*}
  Thus, the algorithm cannot terminate before $t=T-1$.

  Now fix $t\in\{1,\ldots,T-2\}$.  At $t=1$,
  Lemma~\ref{lem:high-level}, the definition of $q(x)$ and the
  assumption that the solver is not lucky ensure that
  $q\left(x^{(0)}\right)=1$. For $t\ge2$,
  Proposition~\ref{prop:forced-level} combined with our assumptions on
  the LP solver gives $q\left(x^{(t)}\right)=t+1$.

  Because $q(x^{(t-1)}) = t$ and $q(x^{(t)}) = t+1$, and both indices
  are at most $T$ for $t \leq T-2$, the supports of $x^{(t-1)}$ and
  $x^{(t)}$ have opposite parity. By
  Proposition~\ref{prop:vertexcases}, because we are switching parity
  the values of the controller variables cannot stay the same across
  iterations: $w_1 \in \{1, \frac{1}{3}, \frac{1}{6}, \dots,
  \frac{1}{3(T+1)}\}$ for odd vertices, and $w_1 \in \{0,
  1-\frac{1}{3}, 1-\frac{1}{6}, \dots, 1-\frac{1}{3(T+1)}\}$ for even
  vertices. Thus, we have
  \begin{equation}
    \label{eq:wdifflb}
    \|w^{(t)} - w^{(t-1)}\|_\infty \geq \frac{1}{3(T+1)}.
  \end{equation}
  Using the IR update rule $w^{(t)} = w^{(t-1)} + \epsilon^t\hat
  w^{(t)}$, we obtain
  \begin{equation*}
    \hat w^{(t)} = \epsilon^{-t} \left( w^{(t)}-w^{(t-1)} \right).
  \end{equation*}
  Combining this with \eqref{eq:wdifflb}, we obtain
  \begin{equation*}
      \|\hat x^{(t)}\|_\infty \geq \left\|\hat w^{(t)}\right\|_\infty 
      = \epsilon^{-t}\left\| w^{(t)}-w^{(t-1)} \right\|_\infty 
      \geq \frac{1}{3(T+1)}\epsilon^{-t} = \bigOm{\epsilon^{-t}}. \qedhere
  \end{equation*}
\end{proof}

\section{Proof of Theorem \ref{thm:main-simplex}}
\label{sec:simplex}
We show that if the limited-precision LP solver in
Algorithm~\ref{alg:ir} is the simplex algorithm using Bland's pivoting
rule, and the solution process of each refining LP is warm-started
with the optimal basis for the previous refining LP, then
Algorithm~\ref{alg:ir} generates a sequence of refining solutions of
exponentially-growing size, thereby proving
Theorem~\ref{thm:main-simplex}. Throughout this section, assume
\begin{equation*}
  T\ge 3,
  \qquad
  0<\epsilon\le \frac12,
  \qquad
  0<\epsfin<\frac34\epsilon^{T-2}.
\end{equation*}

\subsection{Simplex preliminaries and Bland's rule}
\label{sec:bland}
For a primal-feasible basis $B$, let $\bar{c}^B$ be the vector of
reduced costs, and write
\begin{equation*}
  \delta_D(B) := \left\| \min\left\{ \zeroes,\bar{c}^B \right\} \right\|_\infty.
\end{equation*}
The simplex implementation studied in this section determines
optimality with tolerance $\epsilon$ by checking if $\bar c_j^B\ge
-\epsilon$ for every nonbasic variable $j \not\in B$. If this
condition holds, the algorithm stops returning the basic feasible
solution associated with $B$, otherwise it pivots by selecting a
variable with negative reduced cost.

For the reasons stated in Section~\ref{sec:preliminaries}, we can
study the execution of a specific instantiation of the simplex
algorithm on the refining problem by relying on the cumulative problem
\eqref{eq:PT} instead. We state this more precisely below.
\begin{proposition}\label{prop:simplex-equivalence}
  For a given primal-feasible basic pair $(x^{(t-1)},y^{(t-1)})$ for
  \eqref{eq:PT} obtained at iteration $t$ of Algorithm~\ref{alg:ir},
  and a given basis $B$ (not necessarily the one associated with
  $(x^{(t-1)},y^{(t-1)})$), the simplex algorithm applied to the
  refinement problem in standard form at basis $B$ follows the same
  execution path as the simplex algorithm applied to \eqref{eq:PT} at
  basis $B$. If the simplex algorithm uses a nonzero tolerance to
  determine if reduced costs are negative, then the execution path
  remains the same as long as that tolerance is scaled by
  $\epsilon^{-t}$.
\end{proposition}
\begin{proof}
  By Proposition~\ref{prop:cumulative-refining-correspondence} the
  bases of the two problems have a direct correspondence, and the
  vectors of reduced costs for the two problems are related by a
  positive scaling factor $\epsilon^{-t}$. Furthermore, note that in the
  refining problem $\min \{ \hat{c}^{\top} x \mid A x = \hat{b}, x
  \geq \hat{\ell} \}$, at iteration $t$ we have $\hat{b} = \zeroes$
  and $\hat{\ell} = -\epsilon^{-t} x^{(t-1)}$, because
  $(x^{(t-1)},y^{(t-1)})$ is a primal-feasible basic pair for
  \eqref{eq:PT}. If we shift the problem to have zero lower bound and
  put it in standard form with decision variables $\tilde{x} \ge
  \zeroes$, the constraints $A x = \hat{b}$ become
  \begin{equation*}
    A (\tilde{x} + \hat{\ell}) = \zeroes \Longleftrightarrow A \tilde{x} = \epsilon^{-t} A x^{(t-1)} \Longleftrightarrow A \tilde{x} = \epsilon^{-t} b,
  \end{equation*}
  because $x^{(t-1)}$ is primal feasible. Thus, the r.h.s.\ vector of
  the refining problem in standard form is also a positive scaling of
  the original r.h.s.\ vector $b$, with the same scaling factor
  $\epsilon^{-t}$. It follows that as long as we scale the tolerance
  for negative reduced costs by $\epsilon^{-t}$ (see
  Proposition~\ref{prop:cumulative-refining-correspondence} (2.)), the
  simplex algorithm applied to the refining problem and to the
  cumulative problem are completely equivalent.
\end{proof}
Proposition~\ref{prop:simplex-equivalence} lets us study \eqref{eq:PT} directly instead of the refining problem. Next, we recall Bland's pivoting rule.
\begin{definition}[Bland's rule]
  \label{def:blandsrule}
  When more than one variable is a candidate to enter or leave the
  basis, \emph{Bland's rule} prescribes that ties are broken
  as follows:
  \begin{itemize}
  \item Entering variable: among all nonbasic variables with negative reduced cost, choose the one with smallest index.
  \item Leaving variable: among all basic variables that tie at the minimum ratio test, choose the one with the smallest index.
  \end{itemize}
\end{definition}
Let $h_O$ and $h_E$ denote, respectively, the odd and even element of
$\{T,T+1\}$. In this section we assume that the columns of the
standard-form LP \eqref{eq:PT} are permuted according to the
precedence order $\prec$ given by
\begin{equation}
\label{eq:bland-order}
  s_{h_O}\prec s_{h_E}\prec s_1\prec s_2\prec\cdots\prec s_{T-1}
  \prec z_1\prec z_2\prec\cdots\prec z_{T+1}
  \prec w_2\prec w_1.
\end{equation}
Applying Bland's rule with the standard index-based precedence, as
described in Definition~\ref{def:blandsrule}, to the problem with
permuted column order is equivalent to Bland's rule applied to
\eqref{eq:PT} with the precedence order given by $\prec$. Thus, below
for ease of exposition we keep the same variable order as in
\eqref{eq:PT}, but break ties using $\prec$ as indicated in
\eqref{eq:bland-order}.

We assume that every refining LP is warm-started using the final basis
returned by the preceding solve, i.e., the optimal basis of the
previous refining LP. This is a natural choice because this basis is
guaranteed to be primal feasible and it is readily available. Once we
fix the implementation of the simplex algorithm and an initial
solution, we identify a path along the bases and vertices of
\eqref{eq:PT}. For the simplex algorithm with Bland's pivoting rule
and optimality tolerance $\epsilon$ for each refining problem, we
obtain the path that we study in the next subsection.

\subsection{Bases explored following Bland's rule}
\label{sec:blandpath}
For $r\in\{1,\ldots,T-1\}$, let $h(r)$ be the element of $\{T,T+1\}$
whose parity is opposite to that of $r$, and let $g(r)$ be the one
with the same parity.

For each $r\in\{1,\ldots,T-1\}$, we consider the two vertices that
have $z_r > 0, z_j = 0$ for all $j \neq r$. The first vertex, denoted
$u(r) = (w,z,s)$, is the type 2 vertex with $S(z)=\{r\}$, namely
\begin{equation*}
  z_r=\epsilon^r, \qquad z_j=0 \text{ for } j \neq r, \qquad
  w_{p(r)}=\frac13, \qquad w_{3-p(r)}=\frac23.
\end{equation*}
By Lemma~\ref{lem:interior-unique}, $u(r)$ corresponds uniquely to the
canonical basis
\begin{equation*}
  B(u(r)) = \{w_1,w_2,z_r\} \cup
  \{s_j:j\in\{1,\ldots,T+1\}\setminus\{r\}\}.
\end{equation*}
The second vertex, denoted $v(r) = (w,z,s)$, is the type 1 vertex with
$S(z)=\{r\}$, namely
\begin{equation*}
  z_r=\epsilon^r, \qquad z_j=0 \text{ for } j \neq r, \qquad
  w_{p(r)}=1, \qquad w_{3-p(r)}=0.
\end{equation*}
For $v(r)$, Definition~\ref{def:canonicalbasis} identifies the canonical basis
\begin{equation*}
  B(v(r)) =
  \{w_1,w_2,z_r\} \cup \{s_j:j\in\{1,\ldots,T+1\}\setminus\{h(r)\}\}.
\end{equation*}
We introduce an additional symbol for a certain non-canonical basis corresponding to $v(r)$, which will be encountered in the following. For $r\ge2$, we define
\begin{equation*}
  N(v(r)) :=
  \{w_1,w_2,z_{r-1},z_r\} \cup \{s_j:j\in\{1,\ldots,T+1\}\setminus\{r-1,h(r)\}\}.
\end{equation*}
Note that in this basis, $z_{r-1}$ is basic at value zero, whereas $s_{r-1}$ is
nonbasic at value zero.

Next, we study the reduced costs associated with these bases. For
convenience, we write below the instantiation of \eqref{eq:PT} for
$T=5$ and its constraint matrix, which makes the calculations in the
next few proofs easier to verify.

\begin{equation*}
\renewcommand{\arraystretch}{1.5}
\setlength{\arraycolsep}{5pt}
\begin{array}{c|cc|cccccc|cccccc|c}
& w_1 & w_2
& z_1 & z_2 & z_3 & z_4 & z_5 & z_6
& s_1 & s_2 & s_3 & s_4 & s_5 & s_6
& \text{RHS} \\
\hline
\min\; c^\top x
& 0 & 0
& \epsilon^5 & \epsilon^5 & \epsilon^5
& \epsilon^5 & \epsilon^5 & \epsilon^5
& 0 & 0 & 0 & 0 & 0 & 0
& 0 \\[-6pt]
&
& 
& + & + & +
& + & + & +
& & & & & &
& \\[-6pt]
&
& 
& \epsilon^{-1} & \epsilon^{-2} & \epsilon^{-3}
& \epsilon^{-4} & \epsilon^{-5} & \epsilon^{-5}
& & & & & &
& \\
\hline\hline
w_1+w_2=1
& 1 & 1
& 0 & 0 & 0 & 0 & 0 & 0
& 0 & 0 & 0 & 0 & 0 & 0
& 1 \\

\sum_{j=1}^{5}\epsilon^{-j}z_j+\epsilon^{-5}z_6=1
& 0 & 0
& \epsilon^{-1} & \epsilon^{-2} & \epsilon^{-3}
& \epsilon^{-4} & \epsilon^{-5} & \epsilon^{-5}
& 0 & 0 & 0 & 0 & 0 & 0
& 1 \\
\hline
-3w_1+\epsilon^{-1}z_1+s_1=0
& -3 & 0
& \epsilon^{-1} & 0 & 0 & 0 & 0 & 0
& 1 & 0 & 0 & 0 & 0 & 0
& 0 \\

-3w_2+\epsilon^{-2}z_2+s_2=0
& 0 & -3
& 0 & \epsilon^{-2} & 0 & 0 & 0 & 0
& 0 & 1 & 0 & 0 & 0 & 0
& 0 \\

-3w_1+\epsilon^{-3}z_3+s_3=0
& -3 & 0
& 0 & 0 & \epsilon^{-3} & 0 & 0 & 0
& 0 & 0 & 1 & 0 & 0 & 0
& 0 \\

-3w_2+\epsilon^{-4}z_4+s_4=0
& 0 & -3
& 0 & 0 & 0 & \epsilon^{-4} & 0 & 0
& 0 & 0 & 0 & 1 & 0 & 0
& 0 \\

-3w_1+\epsilon^{-5}z_5+s_5=0
& -3 & 0
& 0 & 0 & 0 & 0 & \epsilon^{-5} & 0
& 0 & 0 & 0 & 0 & 1 & 0
& 0 \\

-3w_2+\epsilon^{-5}z_6+s_6=0
& 0 & -3
& 0 & 0 & 0 & 0 & 0 & \epsilon^{-5}
& 0 & 0 & 0 & 0 & 0 & 1
& 0
\end{array}
\end{equation*}

\begin{lemma}[Reduced costs for $B(u(r))$]
  \label{lem:reducedcostur}
  Assume $0<\epsilon\le\frac12$, and let $r \in \{1,\dots,T-1\}$. At
  the basis $B(u(r))$, the reduced costs for nonbasic variables
  satisfy
  \begin{equation*}
    \bar c_{z_j}^{B(u(r))} = \epsilon^T\left(1-\epsilon^{r-\min\{j,T\}}\right)
    \text{ for } j\neq r, \qquad \bar c_{s_r}^{B(u(r))}=0.
  \end{equation*}
\end{lemma}
\begin{proof}
  The dual vector associated with $B(u(r))$ is defined by the linear
  system $A_{B(u(r))}^\top\pi=c_{B(u(r))}$; recall that this system
  has $T+3$ columns. Every slack except $s_r$ is basic, so
  $\pi_{j+2}=0$ for $j \ne r$. Using this fact, the first two rows of
  the system (obtained from the columns corresponding to $w_1, w_2$)
  imply $\pi_1=\pi_{r+2}=0$. The row obtained from the column of $z_r$
  then gives $\pi_2=1+\epsilon^{T+r}$. Applying the formula for
  reduced cost concludes the proof: for each $z_j$, $j \neq r$, we have
  \begin{equation*}
    \bar c_{z_j}^{B(u(r))} = c_{z_j} - \pi^{\top} A_{z_j} = (\epsilon^T + \epsilon^{-\min\{j,T\}}) - (1+\epsilon^{T+r}) \epsilon^{-\min\{j,T\}} = \epsilon^T\left(1-\epsilon^{r-\min\{j,T\}}\right),
  \end{equation*}
  and for $s_r$ we have $c_{s_r}^{B(u(r))}= 0 - 0$.  
\end{proof}
\begin{lemma}[Reduced costs for $N(v(r))$]
  \label{lem:reducedcostvr}
  Assume $0<\epsilon\le\frac12$, and let $r \in \{2,\dots,T-1\}$. At
  the basis $N(v(r))$, the reduced costs for nonbasic variables
  satisfy
  \begin{align*}
    \bar c_{z_j}^{N(v(r))} &= \epsilon^T\left(1-\epsilon^{r-\min\{j,T\}}\right) \text{ for } j\notin\{r-1,r,h(r)\}, \qquad
    \bar c_{z_{h(r)}}^{N(v(r))} = \epsilon^T+\epsilon^{r-1}-2\epsilon^r,\\
    \bar c_{s_{r-1}}^{N(v(r))} &=-(1-\epsilon)\epsilon^{T+r-1}, \qquad
    \bar c_{s_{h(r)}}^{N(v(r))} =(1-\epsilon)\epsilon^{T+r-1}.
  \end{align*}
\end{lemma}
\begin{proof}
  The proof follows the same structure as for
  Lemma~\ref{lem:reducedcostur}. Every slack except $s_{r-1}$ and
  $s_{h(r)}$ is basic, so $\pi_{j+2}=0$ for $j \neq r-1, j \neq
  h(r)$. The remaining part of the system defining $\pi$ gives
  \begin{equation*}
    \pi_2=1+\epsilon^{T+r},
    \qquad
    \pi_{r+1}=(1-\epsilon)\epsilon^{T+r-1},
    \qquad
    \pi_{h(r)+2}=-\pi_{r+1},
    \qquad
    \pi_1=0.
  \end{equation*}
  The reduced cost of $s_j$ is $-\pi_{j+2}$; this proves the formula
  for $\bar c_{s_{r-1}}^{N(v(r))}$ and $\bar c_{h(r)}^{N(v(r))}$. For
  the $z$ variables with $j \notin\{r-1,r,h(r)\}$, we have
  \begin{equation*}
    \bar c_{z_j}^{N(v(r))} = c_{z_j} - \pi^{\top} A_{z_j} = (\epsilon^T + \epsilon^{-\min\{j,T\}}) - (1+\epsilon^{T+r}) \epsilon^{-\min\{j,T\}} = \epsilon^T\left(1-\epsilon^{r-\min\{j,T\}}\right),
  \end{equation*}
  and for $z_{h(r)}$ we have
  \begin{align*}
    \bar c_{z_{h(r)}}^{N(v(r))} &= c_{z_{h(r)}} - \pi^{\top} A_{z_{h(r)}} = (\epsilon^T + \epsilon^{-T}) - (1+\epsilon^{T+r}) \epsilon^{-T} + (1-\epsilon)\epsilon^{T+r-1} \epsilon^{-T} \\
    &= \epsilon^T+\epsilon^{r-1}-2\epsilon^r.
  \end{align*}
\end{proof}
With these two lemmas, we can characterize the dual infeasibility of
the bases of interest.
\begin{corollary}
  \label{cor:dualinfeasr}
  For $r \in \{1,\dots,T-1\}$ we have $\delta_D(B(u(r))) = \epsilon^r
  - \epsilon^T$, and for $r \in \{2,\dots,T-1\}$ we have
  $\delta_D(N(v(r))) = \epsilon^r - \epsilon^T$.
\end{corollary}
\begin{proof}
  Follows immediately from Lemmas \ref{lem:reducedcostur} and
  \ref{lem:reducedcostvr}, noting that the most negative reduced cost
  is for $z_{j}$ with $j = g(r)$, i.e., for the $z$ variable with the
  same parity as $r$ and the largest index (in the case of $N(v(r))$,
  we exclude index $j = h(r)$ as the corresponding reduced cost is
  positive: $\epsilon^T+\epsilon^{r-1}-2\epsilon^r \ge 0$ for
  $\epsilon \le \frac12$).
\end{proof}

We now analyze the basis changes induced by Bland's rule from vertices
$u(r)$ and $v(r)$.
\begin{lemma}[Pivot from $u(r)$]
  \label{lem:urpivot}
  Let $r\in\{1,\ldots,T-2\}$. Suppose we are at iteration $t=r$ of
  Algorithm~\ref{alg:ir}, using the simplex algorithm with Bland's
  pivoting rule and optimality tolerance $\epsilon$, initialized with
  the feasible basis $B(u(r))$, to solve the refining problem. Then
  the LP solver terminates at basis
  \begin{equation*}
    B(u(r)) \cup \{z_{r+1}\} \setminus \{s_{h(r+1)}\} = N(v(r+1)).
  \end{equation*}
\end{lemma}
\begin{proof}
  By Propositions~\ref{prop:cumulative-refining-correspondence} and
  \ref{prop:simplex-equivalence}, and Corollary~\ref{cor:dualinfeasr},
  the dual infeasibility for the refining problem at iteration $t=r$
  is
  \begin{equation*}
    \epsilon^{-r}\delta_D(B(u(r))) = 
    \epsilon^{-r}\left(\epsilon^r-\epsilon^T\right)
    = 1-\epsilon^{T-r} >\epsilon,
  \end{equation*}
  where the inequality uses $r\le T-2$ and $\epsilon\le\frac12$. Thus,
  basis $B(u(r))$ is not $\epsilon$-optimal, and we need to perform at
  least one pivot. The only nonbasic slack is $s_r$, whose reduced
  cost is zero. By Lemma~\ref{lem:reducedcostur}, the $z$ variable
  with negative reduced cost that has the smallest index in the
  $\prec$ order is $r+1$.  Therefore $z_{r+1}$ enters the basis.

  The edge that starts at $u(r)$ and increases $z_{r+1}$ is defined by
  \begin{equation*}
    \epsilon^{-r}z_r=1-\lambda,
    \qquad
    w_{p(r)}=\frac{1-\lambda}{3},
    \qquad
    w_{3-p(r)}=\frac{2+\lambda}{3},
    \qquad
    0\le\lambda\le1.
  \end{equation*}
  Variables $z_r$ and $w_{p(r)}$ are zero at $\lambda=1$; every basic
  slack whose index has the same parity as $r$ also decrease to zero
  as $\lambda$ goes to $1$. All these variables tie in the
  minimum-ratio test, and the one with the smallest index in the
  $\prec$ order is $s_{h(r+1)}$. Therefore $s_{h(r+1)}$ leaves. It is
  immediate to check that at $\lambda = 1$ the above equations yield
  the point $v(r+1)$, and that the new basis is precisely $N(v(r+1))$.
  
  Using Corollary~\ref{cor:dualinfeasr} again, the dual infeasibility
  for the refining problem at the new basis is
  \begin{equation*}
    \epsilon^{-r}\delta_D(N(v(r+1))) = 
    \epsilon^{-r}\left(\epsilon^{r+1}-\epsilon^T\right)
    = \epsilon-\epsilon^{T-r} <\epsilon.
  \end{equation*}
  Thus, the simplex algorithm terminates here.
\end{proof}

\begin{lemma}[Pivots from $v(r)$]
  \label{lem:vrpivot}
  Let $r\in\{2,\ldots,T-2\}$. Suppose we are at iteration $t=r$ of
  Algorithm~\ref{alg:ir}, using the simplex algorithm with Bland's
  pivoting rule and optimality tolerance $\epsilon$, initialized with
  the feasible basis $N(v(r))$, to solve the refining problem. Then
  the LP solver terminates at basis $B(u(r+1))$ after three pivots,
  two of which are degenerate.
\end{lemma}
\begin{proof}
  By Propositions~\ref{prop:cumulative-refining-correspondence} and
  \ref{prop:simplex-equivalence}, and Corollary~\ref{cor:dualinfeasr},
  the dual infeasibility for the refining problem at iteration $t=r$
  is
  \begin{equation*}
    \epsilon^{-r}\delta_D(N(v(r))) = 
    \epsilon^{-r}\left(\epsilon^r-\epsilon^T\right)
    = 1-\epsilon^{T-r} >\epsilon,
  \end{equation*}
  where the inequality uses $r\le T-2$ and $\epsilon\le\frac12$. Thus,
  basis $N(v(r))$ is not $\epsilon$-optimal, and we need to perform at
  least one pivot.  By Lemma~\ref{lem:reducedcostvr}, $s_{r-1}$ has
  negative reduced cost, whereas the other nonbasic slack $s_{h(r)}$
  has positive reduced cost.  Since the slacks precede the $z$
  variables in \eqref{eq:bland-order}, $s_{r-1}$ is the eligible
  entering variable with the smallest index in the $\prec$ order.

  At basis $N(v(r))$, both $z_{h(r)}$ and $s_{h(r)}$ are nonbasic at
  zero. Thus, the constraint $-3 w_{p(h(r))} +
  \epsilon^{-T} z_{h(r)} + s_{h(r)} = 0$ fixes
  $w_{p(h(r))} = 0$ as $s_{r-1}$ increases along its nonbasic
  direction. Because $r-1$ has the same parity as $h(r)$, the row
  containing $s_{r-1}$ reduces to
  \begin{equation*}
    \epsilon^{-(r-1)}z_{r-1}+s_{r-1}=0.
  \end{equation*}
  This is the only row containing $s_{r-1} > 0$, and because $z_{r-1}$
  is basic with value zero, we have a degenerate pivot where $z_{r-1}$
  leaves the basis. The resulting basis is
  \begin{equation*}
    N(v(r)) \cup \{s_{r-1}\} \setminus \{z_{r-1}\} = B(v(r)).
  \end{equation*}

  At this new basis, Lemma~\ref{lem:boundary-reduced-costs} gives
  \begin{equation*}
    \bar c_{z_j}^{B(v(r))} =
    \epsilon^T\left(1-\epsilon^{r-\min\{j,T\}}\right) \text{ for } j \neq r,
  \end{equation*}
  so $\bar c_{z_{g(r)}}^{B(v(r))}=\epsilon^T-\epsilon^r$.  After
  scaling by $\epsilon^{-r}$ because we are at the $r$-th refinement
  iteration (Proposition~\ref{prop:simplex-equivalence}), this reduced
  cost is $\epsilon^{T-r}-1<-\epsilon$, and the basis is not
  $\epsilon$-optimal. The only nonbasic slack is $s_{h(r)}$, whose
  reduced cost is zero, and the $z$ variable with negative reduced
  cost and smallest index in the $\prec$ order is $z_{r+1}$. Hence
  $z_{r+1}$ enters the basis.

  Because $r\le T-2$, $r+1\le T-1$. With an argument similar to the
  one used for the basis $N(v(r))$, we can observe that the nonbasic
  pair $z_{h(r)}=s_{h(r)}=0$ fixes $w_{p(h(r))} = 0$. Since $p(r+1) =
  p(h(r))$, the constraint $-3w_{p(r+1)} + \epsilon^{-(r+1)}z_{r+1} +
  s_{r+1} = 0$ shows that the zero-valued basic slack $s_{r+1}$
  decreases as $z_{r+1}$ increases. Because $s_{r+1} \ge 0$, this is
  another degenerate pivot, where $s_{r+1}$ leaves. The resulting
  basis is
  \begin{equation*}
    B(v(r)) \cup \{z_{r+1}\} \setminus \{s_{r+1}\} = W(v(r)) :=
    \{w_1,w_2,z_r,z_{r+1}\} \cup
    \{s_j:j\in\{1,\ldots,T+1\}\setminus\{r+1,h(r)\}\}.
  \end{equation*}

  At $W(v(r))$, we compute the vector of dual variables as in the proof of Lemma~\ref{lem:reducedcostvr}. With similar calculations, we obtain
  \begin{equation*}
    \pi_2=1+\epsilon^{T+r},
    \qquad
    \pi_{r+3}=-(1-\epsilon)\epsilon^{T+r},
    \qquad
    \pi_{h(r)+2}=(1-\epsilon)\epsilon^{T+r},
  \end{equation*}
  and all remaining components of $\pi$ are zero. Thus,
  \begin{equation*}
    \bar c_{s_{h(r)}}^{W(v(r))} = -(1-\epsilon)\epsilon^{T+r}<0,
    \qquad
    \bar c_{s_{r+1}}^{W(v(r))} = (1-\epsilon)\epsilon^{T+r}>0.
  \end{equation*}
  Because $\pi_{g(r)} = 0$, we have
  \begin{equation*}
    \bar c_{z_{g(r)}}^{W(v(r))}=\epsilon^T-\epsilon^r
  \end{equation*}
  as usual. This implies that $\epsilon^{-r}\delta_D(W(v(r))) =
  1-\epsilon^{T-r} >\epsilon$, so this basis is not
  $\epsilon$-optimal. The variable $s_{h(r)}$ has negative reduced
  cost and it has the lowest index in the order $\prec$, so it enters
  the basis.
  
  Let $s_{h(r)} = \lambda$. Along the edge defined by $s_{h(r)}$, the
  row $h(r)+2$ and the row indexed by $r+3$ (in other words, the rows
  that upper bound $z_{h(r)}$ and $z_{r+1}$, respectively) give
  \begin{equation*}
    w_{p(r+1)}=\frac{\lambda}{3},
    \qquad
    \epsilon^{-(r+1)}z_{r+1}=\lambda.
  \end{equation*}
  The first two constraints then imply
  \begin{equation*}
    w_{p(r)}=1-\frac{\lambda}{3},
    \qquad
    \epsilon^{-r}z_r=1-\lambda.
  \end{equation*}
  For $0 \le \lambda \le 1$, every other basic variable remains
  nonnegative. At $\lambda=1$, $z_r$ is the unique basic variable that
  reaches zero, so it leaves the basis. The new basis is
  \begin{equation*}
    W(v(r)) \cup \{s_{h(r)}\} \setminus \{z_r\} = B(u(r+1)),
  \end{equation*}
  and the associated vertex is precisely $u(r+1)$. By
  Corollary~\ref{cor:dualinfeasr}, the dual infeasibility for this
  basis in the refining problem is
  \begin{equation*}
    \epsilon^{-r}\left(\epsilon^{r+1}-\epsilon^T\right)
    = \epsilon-\epsilon^{T-r} <\epsilon.
  \end{equation*}
  Thus, the simplex algorithm terminates here.
\end{proof}

\subsection{Putting it together}
We now complete the proof of Theorem~\ref{thm:main-simplex}. We
arbitrarily choose the basis returned by the first solve (on
line~\ref{line:cumulprob} of Algorithm~\ref{alg:ir}) to be
\begin{equation}
  \label{eq:initial-simplex-basis}
  B(u(1)) = \{w_1,w_2,z_1\} \cup \{s_j:j\in\{2,\ldots,T+1\}\},
\end{equation}
which corresponds to the vertex of type 2O with $z_1 = \epsilon$; we
check below that this solution satisfies the required tolerances, so
it is a valid solution for the initial solve.

\begin{proof}[Proof of Theorem~\ref{thm:main-simplex}]
  At the initial basis $B(u(1))$, Corollary~\ref{cor:dualinfeasr}
  gives
  \begin{equation*}
    \delta_D(B(u(1)))=\epsilon-\epsilon^T<\epsilon.
  \end{equation*}
  Thus, this basis satisfies the tolerances for the LP solver to
  terminate, returning this basis and the primal point $x^{(0)}=u(1)$.
  Because $T\ge3$, $\epsilon\le\frac12$, and
  $\epsfin<\frac34\epsilon^{T-2}$, we also see that
  \begin{equation*}
    \delta_D^{(0)} = \epsilon-\epsilon^T
    \ge \epsilon(1-\epsilon^2)
    \ge \frac34\epsilon
    \ge \frac34\epsilon^{T-2}
    >
    \epsfin.
  \end{equation*}
  Therefore Algorithm~\ref{alg:ir} proceeds to the first refinement iteration.

  We now argue inductively that, for every $t\in\{1,\ldots,T-2\}$, the
  primal solution and basis returned by the chosen LP solver applied
  to the refining problem are
  \begin{equation}
    (x^{(t)},B^{(t)})
    =
    \begin{cases}
      (v(t+1), N(v(t+1)), & t\textnormal{ odd},\\
      (u(t+1), B(u(t+1)), & t\textnormal{ even}.
    \end{cases}
  \end{equation}
  For odd $t$, the preceding index $t-1$ is even, so the induction
  hypothesis gives $(x^{(t-1)},B^{(t-1)})=(u(t),B(u(t))$.  Applying
  Lemma~\ref{lem:urpivot} with $r=t$ proves the inductive statement
  for $t$ odd.  For even $t$, the induction hypothesis gives
  $(x^{(t-1)},B^{(t-1)})=(v(t), N(v(t)))$, and Lemma~\ref{lem:vrpivot}
  concludes the proof of the inductive statement.

  Along this sequence of solutions, Corollary~\ref{cor:dualinfeasr}
  shows that
  \begin{equation*}
    \delta_D^{(t)}=\epsilon^{t+1}-\epsilon^T.
  \end{equation*}
  Thus, for every $t\in\{1,\ldots,T-3\}$,
  \begin{equation*}
    \delta_D^{(t)} \ge \epsilon^{T-2}-\epsilon^T =
    \epsilon^{T-2}(1-\epsilon^2) \ge \frac34\epsilon^{T-2} > \epsfin.
  \end{equation*}
  This implies that Algorithm~\ref{alg:ir} cannot terminate before
  executing refinement iteration $T-2$.
    
  Finally, $x^{(0)}=u(1)$ and the inductive proof above show
  that the values of the $w$ variables alternate according to
  \begin{equation*}
    (w_1^{(t)},w_2^{(t)})
    =
    \begin{cases}
      \left(\frac13,\frac23\right), & t\textnormal{ even},\\[1mm]
      (0,1), & t\textnormal{ odd},
    \end{cases}
    \qquad
    t\in\{0,\ldots,T-2\}.
  \end{equation*}
  This implies that, for every $t\in\{1,\ldots,T-2\}$,
  \begin{equation*}
    \left\|w^{(t)}-w^{(t-1)}\right\|_\infty=\frac13.
  \end{equation*}
  Using the iterative refinement update
  $w^{(t)}=w^{(t-1)}+\epsilon^t\hat w^{(t)}$, we obtain
  \begin{equation*}
    \|\hat w^{(t)}\|_\infty = \frac13\epsilon^{-t}.
  \end{equation*}
  Therefore
  \begin{equation*}
    \|\hat x^{(t)}\|_\infty \ge \|\hat w^{(t)}\|_\infty =
    \frac13\epsilon^{-t} = \bigOm{\epsilon^{-t}},
  \end{equation*}
  which concludes the proof.
\end{proof}

\section*{Acknowledgments}
The authors are partially supported by ONR award \#N000142312585 and
NSF award \#2430617.

\bibliographystyle{alpha}
\bibliography{sample}

\end{document}